\def\addforarchive{\begin{picture}(0,0)
                   \put(180,132){
                     \put(0,45){{\small\sf math.CT/0512076}}
                     \put(0,30){{\small\sf KCL-MTH-05-15}}
                     \put(0,15){{\small\sf ZMP-HH/05-23}}
                     \put(0, 0){{\small\sf Hamburger$\;$Beitr\"age}}
                     \put(5,-10){{\small\sf zur$\;$Mathematik$\;$Nr.$\;$225}}
                     }
                   \end{picture}} 

\documentclass{conm-p-l}
\usepackage{graphicx,amssymb,bbm}
\usepackage[mathscr]{eucal}

\copyrightinfo{2006}{Ribbon Graphs United}

\newtheorem{theorem}{Theorem}

\newtheorem{proposition}[theorem]{Proposition}
\theoremstyle{definition}

\newtheorem{problem}[theorem]{Problem}
\theoremstyle{remark}
\newtheorem{remark}[theorem]{Remark}

\newcommand\void[1]{}

\numberwithin{equation}{section}

\def\atimes        {\mbox{\small$\#$}}
\def\be            {\begin{equation}}
\def\ee            {\end{equation}}

\def\calc          {{\mathcal C}}

\def\calv          {{\V}}
\def\Cf            {\mbox{\sl Cor}}
\def\CfA           {\ensuremath{\Cf_{\!A}}}
\def\cir           {\,{\circ}\,}
\def\cob           {\mbox{\sl3}cob}
\def\complex       {\ensuremath{\mathbbm C}}
\def\dTFT          {3-d\,TFT}
\def\End           {{\rm End}}
\def\eps           {{\varepsilon}}
\def\eq            {\,{=}\,}
\def\findim        {fi\-ni\-te-di\-men\-si\-o\-nal}
\def\id            {\mbox{\sl id}}
\def\II            {{\mathcal I}}
\def\iN            {\,{\in}\,}
\def\Hom           {{\rm Hom}}
\def\HomAA         {\ensuremath{\mathrm{Hom}_{\!A|A}}}
\def\koerper       {{\Bbbk}} 
\def\M             {\ensuremath{\mathrm M}}
\def\nxt           {\raisebox{.08em}{\rule{.44em}{.44em}}\hspace{.4em}}
\def\obj           {{\mathcal O}bj}
\def\one           {{\bf1}}
\def\oti           {\,{\otimes}\,}

\def\reals         {{\mathbb R}}

\def\Rep           {{\mathcal R}ep}
\def\rmA           {{\rm A}}

\def\rmS           {{\rm S}}
\def\simm          {\cong}
\def\sse           {\scriptsize}
\def\tftc          {\ensuremath{\mbox{\sl tft}_{{\mathcal C}}}}
\def\ti            {\,{\times}\,}
\def\tic           {\,{\otimes_{\scriptscriptstyle\complex}}\,}
\def\To            {\,{\to}\,}
\def\V             {\ensuremath{\mathscr V}}
\def\Vect          {{\mathcal V}ect}
\def\Vectkoerper   {\Vect_\koerper}
\def\X             {{\ensuremath{\mathrm X}}}
\def\Xh            {\ensuremath{\widehat{\mathrm X}}}
\def\zet           {{\mathbb Z}}
\def\Zu            {{\rm Z}}

\begin{document}

\title[Topological and conformal field theory as Frobenius algebras]
{Topological and conformal field theory\\ as Frobenius algebras}

\author{Ingo Runkel}
\address{Dept.\ of Mathematics, King's College London, Strand,
GB--London WC2R 2LS}
\email{ingo@mth.kcl.ac.uk}

\author{Jens Fjelstad}
\address{Cardiff School of Math., Cardiff Univ.,
Senghennydd Road, GB--Cardiff CF24 4AG}
\email{fjelstadj@cardiff.ac.uk}

\author{J\"urgen Fuchs}
\address{Institutionen f\"or fysik, Karlstads Universitet,
Universitetsg.~5, S--65188 Karlstad}
\email{jfuchs@fuchs.tekn.kau.se}

\author{Christoph Schweigert}
\address{Fachbereich Mathematik, Universit\"at Hamburg,
Bundesstra{\rm\ss{}}e 55, D--20146 Hamburg}
\email{schweigert@math.uni-hamburg.de}

\thanks{J.Fj.\ is supported by the EU RTN
network Quantum Spaces-Noncommutative Geometry,
J.Fu.\ by VR under project no.\ 621--2003--2385,
and C.S.\ by the DFG project SCHW 1162/1-1.
}


\begin{abstract}
\addforarchive
Two-dimensional conformal field theory (CFT) can be defined through its
correlation functions. These must satisfy certain consistency conditions
which arise from the cutting of world sheets
along circles or intervals. The construction of a (rational) CFT can be divided
into two steps, of which one is complex-analytic and one purely algebraic.
We realise the algebraic part of the construction with the help of
three-dimensional topological field theory and show that any symmetric
special Frobenius algebra in the appropriate braided monoidal category
gives rise to a solution. A special class of examples is provided by
two-dimensional topological field theories, for which the relevant
monoidal category is the category of vector spaces.  
\end{abstract}

\maketitle

\section{Introduction}

It has been known for some time \cite{dijk,abra,kock} that two-dimensional
{\em topological\/} field theories are the same as finite-dimensional
commutative Frobenius algebras over a field $\koerper$. This correspondence 
can be extended to so-called open/closed topological field theories 
\cite{laza,moor10,sega13,laud,lapf}; the Frobenius algebra
is then no longer, in general, commutative.
It has also been known for some time \cite{fuhk} that a subclass of such
theories, so called lattice topological field theories, can
be constructed from a separable non-commutative Frobenius algebra $A$. In
this case, the commutative Frobenius algebra just mentioned is the centre
of $A$, and only the Morita class of $A$ matters.
Intriguingly, an analogous relationship holds in the much more
involved situation of two-dimensional {\em conformal\/} field theory
\cite{frs1}: the relevant structure is now (a Morita class of) a 
symmetric special Frobenius algebra in a braided monoidal category.
It is this latter relation that is further investigated in this article.

The appearance of Frobenius algebras in two-dimensional conformal
field theory (in {\em CFT\/}, for short) is based on the remarkable fact 
that the problem of constructing a CFT can be separated into a 
``complex-analytic'' and an ``algebraic'' part. Here we are mainly concerned 
with the algebraic part. We formulate it (in section \ref{sec:algebra}) as 
Problem \ref{prob:lem}, and show in Theorem \ref{thm:solve} that indeed
a symmetric special Frobenius algebra provides a solution.

The purpose of section \ref{sec:motivate} is to motivate
Problem \ref{prob:lem} and to explain how it relates to other aspects of
CFT. Accordingly, sections \ref{sec:what}\,--\,\ref{sec:open-closed} contain
a brief outline of CFT. The relation between the complex-analytic
and the algebraic aspects relies on a number of physical ideas which so far
have only partially been cast into a precise mathematical language.
If, however, one accepts those ideas, then solving Problem \ref{prob:lem} 
amounts to the construction of a CFT
with a prescribed rational conformal vertex algebra as its chiral symmetry.

\smallskip

Two-dimensional conformal field theory has its origins in several areas of
physics. In all of these it arises as an effective theory, that is, it applies
to experiments once a suitable limit is taken, like the limit of large system
size or of low energy. For example, the long range behaviour of two-dimensional 
statistical systems at equilibrium 
\cite{Card} is described by a (full, local)
CFT, and the properties of edge states in quantum Hall systems (see e.g\ 
\cite{fpsw}) by chiral conformal field theories.

CFT is also a basic ingredient in string theory, which is a candidate for a
fundamental theory of matter and gravity \cite{POlc}. However, string theory
is not yet sufficiently well understood to unambiguously reproduce or
contradict known experimental results. In this sense the only established
appearance of CFT in physics to date is via effective theories.
Nonetheless the strongest incentive for the mathematical development
of CFT came from its applications to string theory.


\section{Two-dimensional conformal field theory}
\label{sec:motivate}

\subsection{What is a two-dimensional CFT?\hspace*{-4.5pt}}\label{sec:what}

There are several versions of two-dimensional CFT, such as formulations
in terms of nets of von Neumann algebras on two-dimensional Minkowski space
(see e.g.\ \cite{muge7,rehr22}), or in terms of correlation functions in the
euclidean plane (see e.g.\ \cite{gabe8}). The formulation of CFT that will be 
reviewed below has been developed to fit the needs of string theory 
\cite{frsh2,vafa0}. It has been cast in a more axiomatic language in 
\cite{sega8,sega8b}; some recent treatments are \cite{gawe2,Gw2,hukr2}.

In this setting, a CFT is based on two pieces of data: First, a $\zet_+$-graded 
complex vector space $H$ with finite-dimensional homogeneous subspaces $H^{(n)}$
(``the $n$th energy eigenspace in $H$''). $H$ is called the 
{\em space of states\/} of the CFT. By $H^\vee$ we denote the graded dual of 
$H$.  Assuming the existence of a discrete grading is actually a simplification;
there are examples of CFTs which lack such a grading, most importantly
Liouville theory \cite{tesc12}. A more appropriate name for the class of 
theories we study here might thus be {\em compact\/} CFTs.
Also, in the light of the open/closed CFT that will be considered
in section \ref{sec:open-closed} below, what we are discussing is
full (as opposed to chiral, see section \ref{sec:chiral-full} below)
{\em closed\/} CFT, and accordingly $H$ is the space of closed CFT states.

\begin{figure}[t]
  \begin{picture}(300,105)(0,7)
  \put(-10,100){a)}
  \put(0,45){$\X~ = $}
  \put(25,7){
     \put(0,0){\scalebox{0.7}{\includegraphics{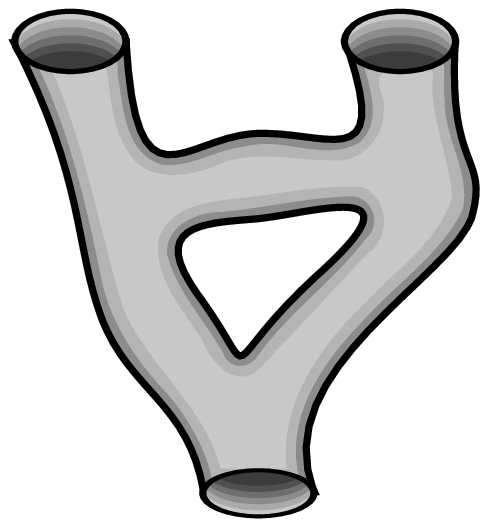}}}
     \put(5,106){\sse out${_1}$}
     \put(72,106){\sse out${_2}$}
     \put(45,-7){\sse in${_1}$}
     }
  \put(170,0){
  \put(0,100){b)}
  \put(0,45){${\rm cut}_f(\X) ~= $}
  \put(45,6){
     \put(0,0){\scalebox{0.7}{\includegraphics{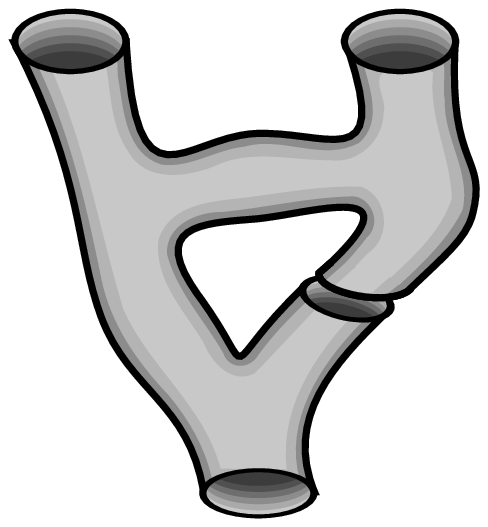}}}
     \put(7,107){\sse out${_1}$}
     \put(75,107){\sse out${_2}$}
     \put(77,36){\sse out${_3}$}
     \put(86,44){\sse in${_2}$}
     \put(45,-7){\sse in${_1}$}
     }}
  \end{picture}
\caption{a) An example of a world sheet and b) the result of
cutting that world sheet along a circle.}
\label{fig:closed-example}
\end{figure}

To describe the second piece of data we need the notion
of a {\em world sheet\/}. A world sheet \X\ is a smooth, compact
two-manifold with parametrised and labelled boundary components
and a Riemannian metric; in this paper we also require that \X\ is oriented.
An example is shown in figure \ref{fig:closed-example}.
The components of the boundary $\partial\X$ are either `incoming' or `outgoing'.
We number the components in the two subclasses separately, and denote by
$|{\rm in}|$ and $|{\rm out}|$ the number of incoming and outgoing components
of $\partial\X$, respectively. Further, every incoming boundary component is
parametrised by an embedding $f\colon~ \rmA^+_\eps \To\X$ of a small annulus
$\rmA^+_\eps \eq \{ z \iN \complex \,|\, 1\,{\le}\,|z|\,{<}\,1{+}\epsilon \}$,
which is is required to be conformal and orientation preserving.
Similarly, outgoing boundary components are parametrised
by embeddings $f\colon~ \rmA^-_\eps \To\X$ of a small annulus
$\rmA^-_\eps \eq \{ z \iN \complex \,|\, 1{-}\eps \,{<}\, |z| \,{\le}\, 1 \}$.
Two world sheets \X\ and $\X'$ are {\em isomorphic\/} iff there exists an 
isometry which is compatible with the boundary labelling and parametrisations.

The second piece of data is a mapping $\Cf$ that assigns to an isomorphism
class $[\X]$ of world sheets an element of the vector space
$\big(H^{|{\rm in}|} \ti (H^\vee)^{|{\rm out}|}\big)^{\!*}$ of multilinear
maps $H \,{\times} \cdots {\times}\, H^\vee \To \complex$.\ For instance, the 
world sheet \X\ in figure \ref{fig:closed-example}\,a) results in a multilinear
map
\be
  \Cf(\X) :\quad  H \ti H^\vee \ti H^\vee \longrightarrow \complex \,.
\ee
$\Cf(\X)$ is called the amplitude, or correlation function, or {\em correlator\/}
for \X.

The assignment $\Cf$ has to fulfil a number of consistency requirements.
To formulate them we need the notion of {\em cutting a world sheet
along a curve\/}. More precisely, let \X\ be a world sheet,
$\rmA_\eps \eq \{ z \iN \complex \,|\, 1{-}\eps\,{<}\,|z|\,{<}\,1{+}\eps \}$
a small annulus, and $f\colon~ \rmA_\eps \To\X$ a conformal orientation
preserving embedding. Then the world sheet ${\rm cut}_f(\X)$ is defined by
cutting along the image of the unit circle, which results in one additional
incoming and one additional outgoing boundary component. For example, from the
world sheet \X\ in figure \ref{fig:closed-example}\,a) one can obtain the world 
sheet ${\rm cut}_f(\X)$ shown in figure \ref{fig:closed-example}\,b).
Thus the correlator of ${\rm cut}_f(\X)$ is a multilinear map
\be
  \Cf({\rm cut}_f(\X)) :\quad
  H \ti H \ti H^\vee \ti H^\vee \ti H^\vee \longrightarrow \complex \,.
\ee
We also define a partial evaluation (or `trace') operation ${\rm tr}_{\rm last}$
from $\big(H^n \ti (H^\vee)^m)\big)^{\!*}$ to $\big(H^{n-1} \ti (H^\vee)^{m-1})
\big)^{\!*}$ as evaluation of the last out- on the last in-com\-po\-nent.\,%
  \footnote{~For doing so, one first performs the evaluation up to some fixed 
  grade $N$ and then takes the limit $N\To\infty$.
  Property (C1) implies that for $\Cf({\rm cut}_f(\X))$ this limit exists.
  Since the embedding $f{:}~ \rmA_\eps \To \X$ does not appear on the left
  hand side of the equality in (C1), it also follows that
  the trace is independent of the parametrisation of the embedded annulus.}

The data $H$, $\Cf$ define a closed CFT iff the following
conditions are satisfied:
\begin{itemize}
\item[{\bf (C1)}] {\sl Factorisation\/}:
  \\[-14pt]~
$$
  \Cf(\X) = {\rm tr}_{\rm last} \Cf({\rm cut}_f(\X))
$$
 for every world sheet \X\ and every embedding $f$.\\[-.7em]~
\item[{\bf (C2)}] {\sl Weyl transformations\/}:
For any two metrics $g$ and $g'$ on \X\ that are related by
$g'_p \eq {\rm e}^{\sigma(p)} g_p$
with some smooth function $\sigma\colon~ \X \To \reals$, one has
$$
  \Cf(\X\;{\rm with\;metric}\;g) = {\rm e}^{c S[\sigma]}
  \Cf(\X\;{\rm with\;metric}\;g') \,,
$$
where $c \iN \complex$ and $S[\sigma]$ is the Liouville action
(see \cite{gawe2,Gw2} for details).
  \end{itemize}
The number $c$ appearing here is called the {\em central charge\/} of the CFT.


\subsection{Special case: Two-dimensional topological field theory}

A class of examples for closed CFTs is provided by two-dimensional
topological field theories \cite{atiy6}; detailed
expositions can be found in \cite{quin,BAki,kock}.

Let $B$ be a finite-dimensional commutative Frobenius algebra over $\complex$. 
Then a CFT in the sense described above is obtained upon setting $H \eq H^{(0)}
\eq B$ and constructing $\Cf$ by cutting a world sheet into elementary building 
blocks on which $\Cf$ is defined in terms of the defining data of $B$. 
For example, consider the world sheets
$$
  \begin{picture}(270,60)(0,0)
  \put(25,35){$X_m = $}
  \put(60,5){
     \put(0,0){\scalebox{.7}{\includegraphics{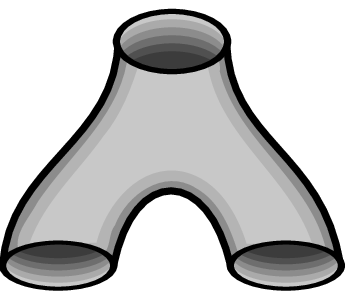}}}
     \put(5,-7){\sse in$_1$}
     \put(54,-7){\sse in$_2$}
     \put(28,60){\sse out$_1$}
     }
  \put(180,35){$X_e = $}
  \put(212,24){
     \put(0,0){\scalebox{.7}{\includegraphics{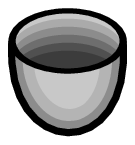}}}
     \put(5,29){\sse out$_1$}
     }
  \end{picture}
$$
The associated correlators $\Cf(\X_m)\colon~ B\ti B\ti B^\vee\To\complex$ and
$\Cf(\X_e)\colon~ B^\vee \To \complex$ are then defined in terms
of the multiplication $m\colon~ B\ti B \To B$ of $B$, and by the
unit $e$ of $B$ (via the map $\eta\colon~ \lambda \,{\mapsto}\, \lambda\,e$
from $\complex$ to $B$), respectively.

The pair $B,\Cf$ obeys condition {\bf(C2)}, with $c\eq0$, for the trivial reason
that $\Cf$ does not depend on the metric at all.
Independence of the world sheet metric also explains the
name {\em topological\/} field theory. That $B$ and $\Cf$ obey {\bf(C1)} follows
from, and is indeed equivalent to, the defining properties of the commutative
Frobenius algebra $B$; this result is due to \cite{dijk} and \cite{abra},
see also \cite{voro2,quin,BAki,kock}.


\subsection{Open/closed CFT}\label{sec:open-closed}

The point of view from which  open/closed CFT is described below is again 
inspired by string theory, namely from models that include both
closed strings, which are circles, and open strings, which are intervals
\cite{POlc,Ansa,sasT2}. In the physics literature, the study of open/closed CFT
reaches back to \cite{card3,cale,lewe3}. The presentation below is in the spirit
of the one given in section \ref{sec:what}; other mathematical approaches along
similar lines have recently been developed in \cite{huan14,huko,hukr3}.  

An open/closed CFT is defined similarly to section \ref{sec:what}, but for a 
wider class of world sheets. Namely, the manifolds are
now allowed to have corners, for example
\be
  \begin{picture}(200,141)(0,0)
  \put(40,5){
     \put(0,0){\scalebox{.7}{\includegraphics{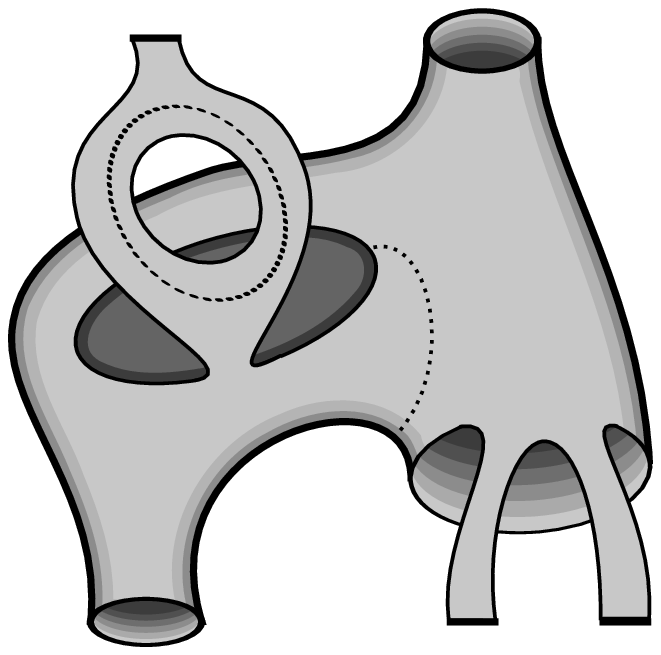}}}
     \put(18,-7){\sse c-in$_1$}
     \put(85,-3){\sse o-in$_1$}
     \put(115,-3){\sse o-in$_2$}
     \put(20,128){\sse o-out$_1$}
     \put(85,132){\sse c-out$_2$}
     }
  \end{picture}
  \label{eq:oc-world}
\ee
In particular there are now three kinds of boundaries, to which we will refer
as {\em closed state boundaries\/}, {\em open state boundaries\/}, and
{\em physical\/} (or {\em external\/}) {\em boundaries\/}, respectively.
Closed state boundaries occurred already in closed CFT in section
\ref{sec:what}. They are circles with a parametrised neighbourhood,
and they are labelled either as incoming or as outgoing and are numbered.
Open state boundaries are intervals that lie between two corners;
they, too, are labelled either as incoming or as outgoing and are numbered.
In addition they are parametrised by a (conformal, orientation
and boundary preserving) embedding of a small
half-annulus $\rmS^+_\eps \eq \rmA^+_\eps \,{\cap}\, \mathbb{H}$ for incoming
boundaries and $\rmS^-_\eps \eq \rmA^-_\eps\,{\cap}\,\mathbb{H}$ for outgoing
boundaries, with $\mathbb{H}$ the upper half plane. For example,
$$
  \begin{picture}(220,59)(0,0)
  \put(0,-14){
    \put(0,0){\scalebox{.6}{\includegraphics{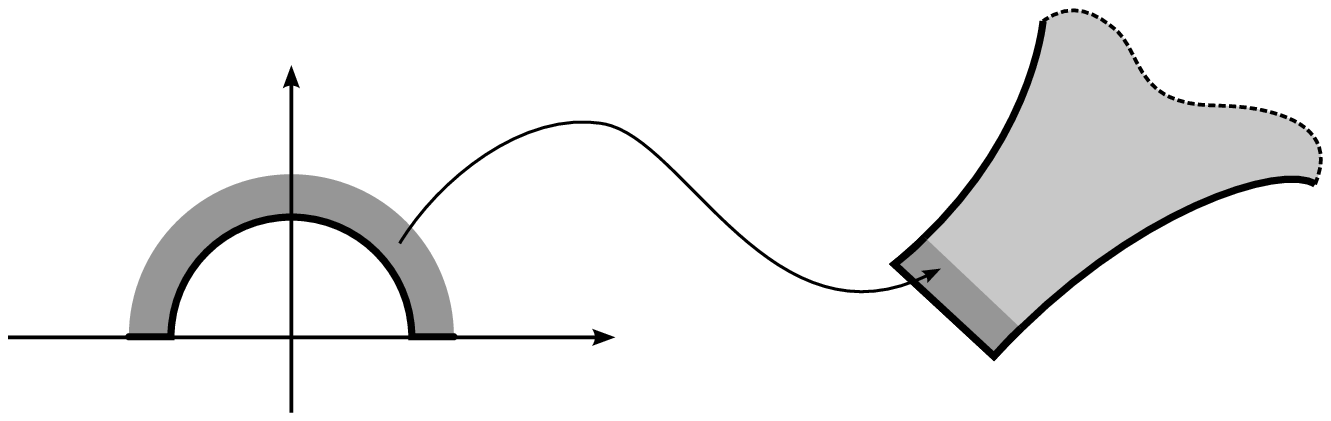}}}
    \put(52,68){\sse Im($z$)}
    \put(88,17){\sse Re($z$)}
    \put(20,47){\sse $S^+_\eps$}
    \put(120,53){\sse $f$}
    \put(146,25){\sse o-in$_k$}
    \put(155,70){\sse phys.}
    \put(155,62){\sse bnd.}
    \put(194,35){\sse phys.\,bnd.}
    }
  \end{picture}
$$
Physical boundaries are not labelled\,\footnote{~%
  Here we make again a simplification: a physical boundary must actually be
  labelled by a so-called {\em boundary condition\/}. To keep the exposition
  short, for simplicity we think of all physical boundaries as having one
  and the same boundary condition, and suppress the label for
  that boundary condition in our notation.}
and not parametrised. They are either circles or intervals between two corners.
At each corner, a physical and an open state boundary meet.
Two open/closed world sheets are isomorphic if they are isometric
in a way compatible with labelling and parametrisation.

The basic data for an open/closed CFT are two graded $\zet_+$-graded
vector spaces $H_{\rm cl}$ and $H_{\rm op}$ (with \findim\ homogeneous subspaces)
and a mapping $\Cf$, called again correlator or amplitude, that assigns to an
isomorphism class of open/\linebreak[0]clo\-sed world sheets a multilinear map.
For example, to the world sheet \eqref{eq:oc-world} a multilinear map
\be
  \Cf(\X) :\quad H_{\rm cl} \ti H_{\rm op} \ti H_{\rm op}
  \ti H_{\rm cl}^\vee \ti H_{\rm op}^\vee \rightarrow \complex 
  \label{eq:corr-as-map}
\ee
is assigned. The conditions for the triple $H_{\rm cl}, H_{\rm op}, \Cf$ to 
define an open/closed CFT are the same as {\bf(C1)} and {\bf(C2)} in section
\ref{sec:what}, except that now there are two distinct ways
of cutting a surface. One can either embed a small annulus
$\rmA_\eps$, in which case cutting along the image of the unit circle results
in two new closed state boundaries, or one can embed a small half-annulus 
$\rmS_\eps \eq \rmA_\eps \,{\cap}\, \mathbb{H}$, so that cutting along the 
image of the unit half-circle results in two new open state boundaries. 
The dashed lines in \eqref{eq:oc-world} show two possible cutting paths.

As there are now two distinct cutting procedures, there are also two partial
trace operations ${\rm tr}_{\rm last,cl}$ and ${\rm tr}_{\rm last,op}$;
the former is the evaluation for the last pair of $H_{\rm cl}$
and $H_{\rm cl}^\vee$, while the latter is the one for the last pair
of $H_{\rm op}$ and $H_{\rm op}^\vee$. Condition {\bf(C1)} is imposed
with respect to each of these two operations.

\medskip

Again, two-dimensional topological field theory provides examples for such a
structure. In this particular class of theories, one considers a 
finite-dimensional symmetric special Frobenius algebra $A$ over \complex, and 
obtains an open/closed CFT with $H_{\rm op} \,{:=}\, A$; the results of 
\cite{frs1} take a particularly simple form in this case, and one finds 
$H_{\rm cl} \,{:=}\, Z(A)$, the centre of $A$. Open/closed 2-d\,TFT has been studied in 
\cite{laza,moor10,sega13} and in the detailed work \cite{lapf}. As before, the 
correlators $\Cf$ are constructed by decomposing world sheets into simple 
building blocks; we refrain from going into any more detail here. Again, the 
open/closed CFT obtained in this way does not depend on the world sheet metric.

As we will see in section \ref{sec:algebra}, by considering a symmetric special 
Frobenius algebra not in the category $\Vect_\complex$ of finite-dimensional 
complex vector spaces, but rather in a more general braided monoidal category, 
we are able to describe also open/closed CFTs that are not topological.


\subsection{Chiral CFT and full CFT}\label{sec:chiral-full}

One way to achieve a construction of open/\linebreak[0]closed CFT is to proceed
in two steps. 
As far as the correlators are concerned, the first step may be stated as
   \begin{itemize}\item[{}]
   ``\,{\sl find all possible functions consistent with a given symmetry\/},\,''
   \end{itemize}
while the second step can be summarised as
   \begin{itemize}\item[{}]\hspace*{-2em}
   ``\,{\sl among all functions with the correct symmetry, find the 
   correlators\/}.\,''
   \end{itemize}
This approach has proven particularly powerful in what is known (see below) as 
{\em rational\/} CFT; it has its roots in \cite{bepz,mose3}.
Let us describe the basic ingredients of this approach.

\smallskip

The symmetries of a CFT can be described by a conformal vertex algebra $\calv$.
To a tuple $(\lambda_1,\lambda_2,...\,,\lambda_m)$ of $\calv$-modules and a 
complex curve $C$ with $m$ marked points one can associate a vector space 
$\mathcal{B}_C(\lambda_1,\lambda_2,...\,,\lambda_m)$ of {\em conformal 
blocks\/}, which is a subspace of the space of multilinear maps from 
$\lambda_1 \ti \lambda_2 \ti \cdots \ti \lambda_m$ to $\complex$. 
This subspace can be constructed explicitly as the space of 
invariants in the algebraic dual of the algebraic tensor product 
$\lambda_1 \tic \lambda_2 \tic \cdots \tic \lambda_m$ with respect to the
dual of the action of $\calv$ on this tensor product.\,\footnote{~%
  In the particular case that $\calv$ is the vertex algebra associated with
  an untwisted affine Lie algebra $\mathfrak g^{\scriptscriptstyle(1)}$, 
  the action of $\calv$ on $\lambda_1\tic\lambda_2\tic\cdots\tic\lambda_m$
  is through the Lie algebra $\mathfrak g\tic\mathcal F$, with
  $\mathfrak g$ the  fi\-ni\-te-di\-men\-si\-o\-nal simple Lie 
  algebra underlying $\mathfrak g^{\scriptscriptstyle(1)}$
  and $\mathcal F$ the associative algebra of meromorphic functions on
  $C$ that have at most finite order poles at the marked points and are
  holomorphic everywhere else \cite{tsuy}. 
  Details for the general case can be found in \cite{FRbe}.}
For a fixed $m$-tuple $(\lambda_1,\lambda_2,...\,,\lambda_m)$ and fixed genus 
$g$ of $C$ the spaces $\mathcal{B}_C(\lambda_1,\lambda_2,...\,,\lambda_m)$ of 
conformal blocks form a vector bundle $\mathcal{B}(\lambda_1,\lambda_2,...\,,
\lambda_m)$ over the moduli space $\mathcal{M}_{g,m}$ of complex curves of 
genus $g$ with $m$ marked points. Each such bundle of conformal blocks is 
equipped with a projectively flat connection.

CFT on complex curves, with the associated system of bundles of conformal 
blocks, is called {\em chiral\/} CFT. Step 1 in the construction of a 
{\em full\/} CFT, i.e.\ of an open/closed CFT in the sense of section 
\ref{sec:open-closed}, consists in determining the corresponding chiral CFT.
The chiral CFT does not fix the correlators uniquely,
but it does determine, for each world sheet $\X$, a subspace
\be
  V(\X) \subset \big( H_{\rm cl} \,{\times} \cdots
  {\times}\, H_{\rm op}^\vee \big)^* 
  \label{eq:cft-V(X)}
\ee
to which the correlator must belong. This subspace $V(\X)$ is given as a space 
of conformal blocks. The relevant surface is, however, not the world sheet \X\ 
itself (which e.g.\ may have a non-empty boundary), but rather its complex 
double \Xh, which is a complex curve, having two marked points for each bulk 
insertion and one marked point for each boundary insertion (see section 
\ref{sec:X-V(X)} for the topological analogue of the complex double). 
$V(\X)$ is the space of conformal blocks on $\Xh$. That bulk insertions 
result in two marked points while each boundary insertion only leads to one 
marked point is in accordance with the fact that one finds an action of 
$\calv\ti\calv$ on $H_{\rm cl}$ and an action of $\calv$ on $H_{\rm op}$.

The symmetries of the CFT, as encoded in the conformal vertex algebra $\calv$,
force the correlator $\Cf(\X)$ of the full CFT to be an element of the space 
$V(\X)$,
\be
  \Cf(\X) \in V(\X) \,.
  \label{eq:corr-in-blocks}
\ee 
In quantum field theoretic terms, the symmetries give rise to
partial differential equations, known as chiral Ward identities,
that must be obeyed by the correlators; they can be formulated as the
condition on the sections of the bundles 
$\mathcal{B}(\lambda_1,\lambda_2,...\,,\lambda_m)$ to be covariantly constant.
   
\smallskip

In step 2 of the procedure one constructs a full CFT
from the chiral CFT obtained in step 1. 
The crucial point is now that for achieving this, everything that is
needed as input from step 1 is already encoded in
the representation category $\Rep(\calv)$ of $\calv$.
Under suitable conditions on $\calv$, $\Rep(\calv)$
is ribbon \cite{hule3.5} and even a modular tensor category \cite{huan21}. 
We refer to a chiral CFT for which $\Rep(\calv)$ is modular
as a {\em rational\/} chiral CFT; it is this class of theories
to which the construction in section \ref{sec:algebra} applies.
Carrying out step 2 can then be formulated entirely as an algebraic
problem in the modular tensor category $\calc$, 
without the need to make any further
reference to the complex-analytic considerations involved in step 1.

That the algebraic problem in step 2 does indeed result in an open/closed CFT
as defined in section \ref{sec:open-closed} relies on several properties of
the spaces $V(\X)$ under cutting of \X\ and under deformation of the
metric (and hence of the complex structure) on \X. These are
natural from the physical perspective, but a proof that they are indeed
fulfilled in any rational CFT is still missing, even though some pertinent
issues have recently been clarified \cite{FRbe,huko2}.
For a more detailed account and references see e.g.\ section 5 of \cite{frs4}.

\smallskip

In section \ref{sec:algebra} we formulate the algebraic problem to be solved 
in the category $\calc$. We do so without making any further reference to the
underlying chiral CFT. Indeed, we do not need to assume that $\calc$ is the
category of representations of a vertex algebra.


\section{The construction of full CFT as an algebraic problem}
\label{sec:algebra}

In this section $\calc$ is a modular tensor category; its defining
properties (which we take to be slightly stronger than in \cite{TUra})
are e.g.\ listed in section 2 of \cite{us}. We take the ground ring 
$\Hom(\one,\one)$ to be an algebraically closed field $\koerper$.

\subsection{Statement of the problem}\label{sec:statement}

We define a {\em topological world sheet\/} \X\ in the same way as a world 
sheet in sections \ref{sec:what} and \ref{sec:open-closed}, except that \X\ 
does not come equipped with a metric. Accordingly, the parametrisations are 
just required to be orientation preserving and continuous, and an isomorphism 
$\X \,{\overset\simm\longrightarrow}\,{\rm Y}$ of topological world sheets is 
a continuous map compatible with orientations and parametrisations.  For $f{:}~ 
\rmA_\eps \To\X$ or $f{:}~\rmS_\eps \To \X$ continuous orientation preserving 
embeddings, the cut topological world sheet ${\rm cut}_f(\X)$ is defined in 
the same way as in sections \ref{sec:what} and \ref{sec:open-closed}, too.
{}From now on we will only consider topological world sheets; for brevity we 
will therefore omit the explicit mentioning of the qualification `topological.'

We denote by $\overline{\calc}$ the dual modular tensor category of $\calc$,
that is, $\calc$ with the brai\-ding replaced by the inverse braiding,
see \cite[section 7]{muge9} and \cite[section 6.2]{C}. Further, let
$\calc \,{\boxtimes}\, \overline{\calc}$ be the product of $\calc$ and
$\overline{\calc}$ in the sense of enriched (over $\Vectkoerper$) category
theory, i.e.\ the modular tensor category obtained by idempotent completion
of the category whose objects are pairs of objects of $\calc$ and
$\overline{\calc}$ and whose morphism spaces are tensor products over
$\koerper$ of the morphism spaces of $\calc$ and $\overline{\calc}$.

With these concepts, and using mappings $\varphi^\sharp$,
${\rm tr}_{\rm last}$ and vector spaces $V(X)$ to be defined below,
the task of constructing a full CFT can be summarised as follows.

\begin{problem} \label{prob:lem}
~\\[.2em]
(i)\hspace*{6pt}Select objects
$H_{\rm cl}\iN\obj(\calc\,{\boxtimes}\,\overline{\calc})$ and
$H_{\rm op}\iN\obj(\calc)$\,.
~\\[.4em]
(ii) Find an assignment $\X \,{\mapsto}\, \Cf(\X) \iN V(\X)$ such that
\begin{itemize}
\item[{\bf (A1)}]
For any isomorphism $\varphi\colon~ \X\,{\overset\simm\longrightarrow}\,{\rm Y}$
one has $\Cf({\rm Y}) \eq \varphi^\sharp \Cf(\X)\,$.
\\[-.7em]~
\item[{\bf (A2)}]
For any embedding $f$ of $\rmA_\eps$ or $\rmS_\eps$ into \X,
$\Cf(\X) \eq {\rm tr}_{\rm last}(\Cf({\rm cut}_f(\X) )\,$.
\end{itemize}
\end{problem}

\noindent
Once we have chosen $H_{\rm cl}$, we further select objects $B_l$ and $B_r$ of 
$\calc$ in such a way that $H_{\rm cl}$ is a subobject of $B_l \,{\times}\, 
\overline{B_r}$. With the help of these objects we can give a prescription, 
to be detailed in sections \ref{sec:X-V(X)} and \ref{sec:trlast} below, how 
to obtain from $\calc, H_{\rm cl}$ and $H_{\rm op}$ the
quantities $\varphi^\sharp$, ${\rm tr}_{\rm last}$ and $V(X)$ entering in the 
formulation of part (ii) of the problem, i.e.\ how to obtain an assignment of
\\[2pt] $\bullet$
a $\koerper$-vector space $V(\X)$ to every world sheet \X;
\\ $\bullet$
a vector space isomorphism
$\varphi^\sharp\colon~ V(\X) \,{\overset\simm\longrightarrow}\, V({\rm Y})$
to every isomorphism $\varphi\colon~ \X \,{\overset\simm\longrightarrow}\,
{\rm Y}$;
\\[2pt] $\bullet$
a map $\,{\rm tr}_{\rm last}\colon~ V({\rm cut}_f(\X)) \To V(\X)$
to every embedding $f$ of $\rmA_\eps$ or $\rmS_\eps$ into \X.

\medskip

\void{
Let us make a few comments on how Problem \ref{prob:lem} relates to the
discussion in section \ref{sec:motivate}. As the notation
suggests, the vector space $V(\X)$ appearing above is the analogue of the
subspace (\ref{eq:cft-V(X)}), and $\Cf(\X)$ plays the role of the correlator.
A solution to Problem \ref{prob:lem} can be translated to elements of
$\Hom( H_{\rm cl} \,{\otimes} \cdots {\otimes}\, H_{\rm op}^\vee,\complex)$;
this yields a solution to conditions {\bf(C1)} and {\bf(C2)} for an open/closed
CFT only if the conformal blocks of the underlying vertex algebras are
sufficiently well-behaved. To find necessary and sufficient conditions on
a vertex algebra for this to be the case remains a challenge.
}

\begin{remark}
Let us make a few comments on how Problem \ref{prob:lem} relates to the
discussion in section \ref{sec:motivate}. As the notation suggests, 
$H_{\rm cl}$ and $H_{\rm op}$ correspond the spaces of states of an 
open/closed CFT. But the space $V(\X)$ in 
Problem \ref{prob:lem} does not correspond to the space of conformal blocks
on the double of $\X$ that appeared in the discussion in section
\ref{sec:chiral-full} (this would require a complex structure),
but rather to the space of flat sections in the relevant vector bundle
of conformal blocks\,\footnote{~%
  More precisely, as the bundle 
  $\mathcal{B}(\lambda_1,\lambda_2,...\,,\lambda_m)$ of conformal blocks 
  over $\mathcal{M}_{g,m}$ is equipped only with a projectively flat connection,
  one must consider flat sections in the projective bundle
  $\mathbb{P}\mathcal{B}(\lambda_1,\lambda_2,...\,,\lambda_m)$.}.
The assignment $\Cf$ in Problem \ref{prob:lem} is the analogue of the correlator
on the CFT side; as in equation \eqref{eq:corr-in-blocks} it is an element in 
the vector space $V(\X)$.

Note that in Problem \ref{prob:lem} we have been careful to {\em not\/} regard
$\Cf(\X)$ as a morphism between tensor products of state spaces, as one
might be tempted to do for having an analogue of \eqref{eq:corr-as-map}. 
(In such an approach the morphisms could only depend on the topology
of \X, and hence one could only describe two-dimensional
{\em topological\/} field theories.) Instead, 
$\Cf(\X)$ is an element of $V(\X)$ and thus, in conjunction with step 1,
is a section of a bundle over $\mathcal{M}_{g,m}$. When endowing
the topological world sheet $\X$ with a complex structure
one selects a point in $\mathcal{M}_{g,m}$ and thereby a
specific conformal block. In this way the dependence of the correlators on 
the complex structure on the world sheet is recovered.
\end{remark}


\subsection{Three-dimensional topological field theory}

To obtain the prescription announced after Problem \ref{prob:lem}, and to
arrive (in section \ref{sec:corr}) at the solution to the Problem, we use as a
crucial tool three-dimensional topological field theory (\dTFT). It
appeared originally in the guise of Chern-Simons field theory 
\cite{schw',witt27,frki2}, and was mathematically developed
in \cite{retu,retu2,tura6,TUra}; 
for reviews see \cite{BAki,KArt} or section 2 of \cite{frs1}.

Let us give a brief outline following \cite{TUra,BAki} to set the notation.
A \dTFT\ furnishes a monoidal functor \tftc\ from a cobordism category 
$\cob_\calc$ to the category $\Vectkoerper$ of finite-dimensional 
$\koerper$-vector spaces. The objects of $\cob_\calc$ are {\em extended 
surfaces\/} $E$, that is, oriented, compact, closed two-manifolds with a 
finite number of disjoint marked arcs and a distinguished Lagrangian subspace 
of $H_1(E,\reals)$. The arcs are marked by pairs $(U,\pm)$ with $U$ an object 
of $\calc$. The following is an example of an extended surface:
$$ \begin{picture}(115,61)(0,0)
    \put(-35,24){$E~~=$}
    \put(0,0){
      \put(0,0){\scalebox{.37}{\includegraphics{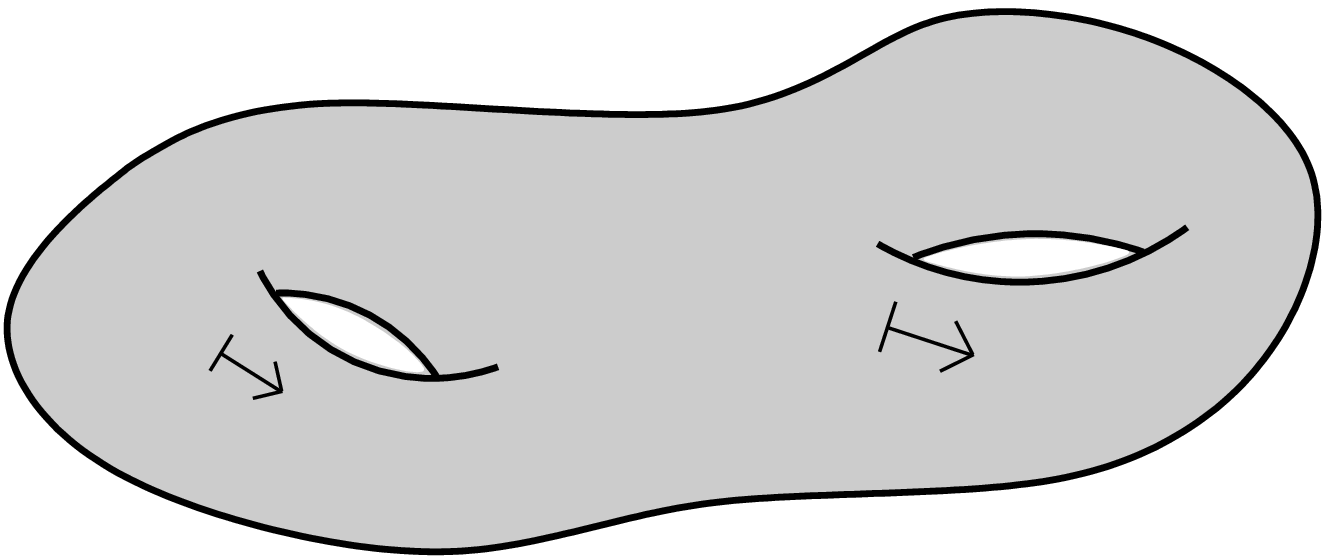}}}
      \put(20,10){\sse $(U,+)$}
      \put(90,15){\sse $(V,-)$}
    }
\end{picture}$$
There is a tensor product on $\cob_\calc$, which on objects
is given by disjoint union; the tensor unit is the empty set.

The morphisms of $\cob_\calc$ are homotopy classes of
{\em extended cobordisms\/} $\M\colon~ E\To E'$. An extended cobordism is a
compact, oriented three-manifolds with boundary given by the disjoint
union of $E$ and $E'$, with an embedded oriented ribbon graph, and with a weight
$m\iN\zet$. The ribbon graph consists of ribbons labelled by objects of
$\calc$ and of coupons labelled by morphisms of $\calc$. The ribbons
end either on coupons or on the marked arcs of $E$ or $E'$. Here is an example:
$$ \begin{picture}(70,168)(30,0)
    \put(-30,87){$M~~=$}
    \put(0,0){
      \put(0,0){\scalebox{.38}{\includegraphics{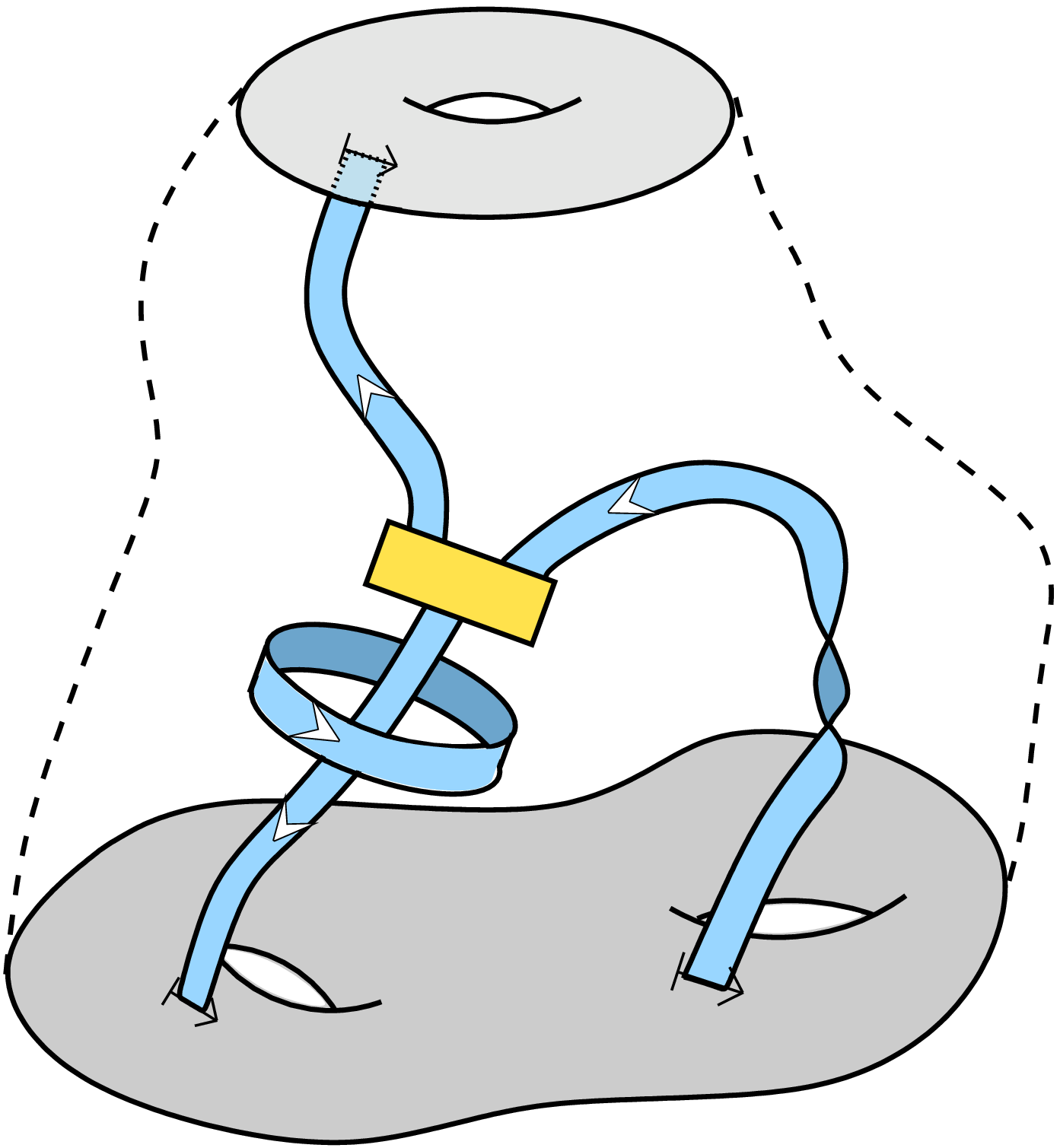}}}
      \put(20,10){\sse $(U,+)$}
      \put(90,15){\sse $(V,-)$}
      \put(60,140){\sse $(W,-)$}
      \put(42,38){\sse $U$}
      \put(90,101){\sse $V$}
      \put(54,115){\sse $W$}
      \put(30,70){\sse $R$}
      \put(63,80){\sse $f$}
      \put(-11,21){$E$}
      \put(22,152){$E'$}
    }
\end{picture}$$
Here $R,U,V,W$ are objects of $\calc$ and $f \iN \Hom(U^\vee, W \oti V^\vee)$.
The identity morphism $\id_E$ is given by $E \,{\times}\, [0,1]$, the
composition of morphisms is just the gluing of the three-manifolds, together
with a rule to compute the new integer weight (see \cite{TUra} for details).
The tensor product of morphisms is again given by disjoint union.


\subsection{The assignment $X \,{\mapsto}\, V(X)$} \label{sec:X-V(X)}

After these preliminaries on \dTFT\ we can present the first part of the
prescription announced after Problem \ref{prob:lem}, namely the assignment
$\X \,{\mapsto}\, V(\X)$.  We first need to select $B_l,\, B_r \iN \obj(\calc)$
such that $H_{\rm cl}$ is a subobject of $B_l \,{\times}\,\overline{B_r}
\iN \obj(\calc \,{\boxtimes}\, \overline\calc)$. This is a technical point:
it is indeed not necessary to introduce these objects, but, as it turns out,
passing to a larger object in $\calc \,{\boxtimes}\, \overline\calc$ that has
a factorised form simplifies the presentation.\,\footnote{~%
    In \cite{fjfrs} a different approach has been followed, in which the object
    $B_l\,{\times}\,\overline{B_r}$ does not appear directly. To recover that
    approach from the present one, one has to decompose
    $B_l\,{\times}\,\overline{B_r}$ into a direct sum of simple subobjects 
    and retain only those which also appear in $H_{\rm cl}$.}
Anticipating the result of Theorem \ref{thm:solve} below, let us also
abbreviate $A \,{\equiv}\, H_{\rm op} \iN \obj(\calc)$.

Consider a world sheet \X\ part of which looks as follows:
$$ \begin{picture}(110,81)(0,0)
    \put(-34,46){$X~~=$}
    \put(0,0){\scalebox{.34}{\includegraphics{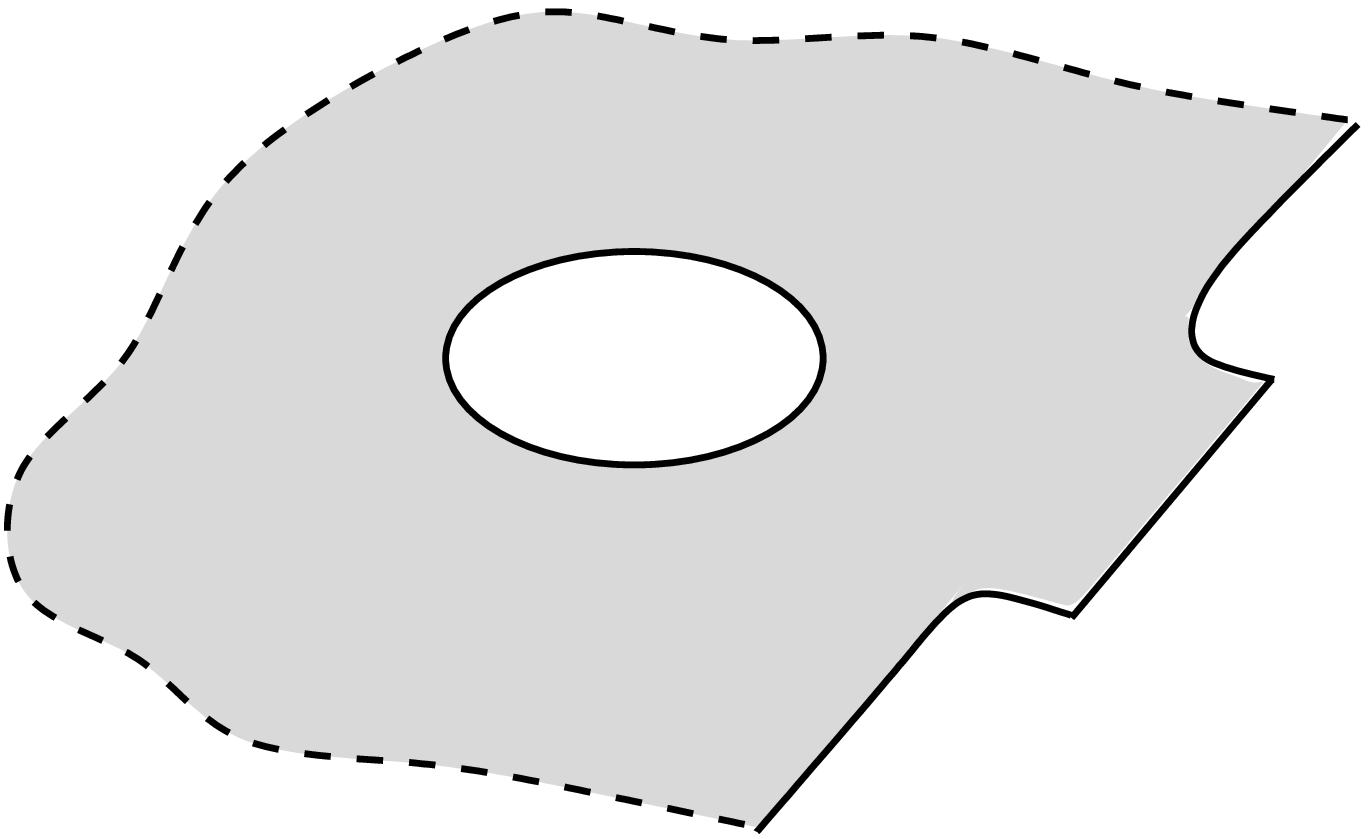}}}
    \put(50,48){\sse c-in$_m$}
    \put(115,30){\sse o-in$_k$}
    \put(87,10){\sse phys.}
    \put(82, 3){\sse bnd.}
    \put(130,62){\sse phys.}
    \put(125,55){\sse bnd.}
\end{picture}$$
The aim is to construct from \X\ an extended surface \Xh, the {\em double\/}
of \X. This will be done in two steps. First, using the parametrisation of the
state boundaries, we glue a disc with an embedded arc to each closed state
boundary, and a semi-disc with an embedded arc on the boundary to each open
state boundary. This results in a surface $\overline\X$ with only physical
boundaries. In the example above,
$$ \begin{picture}(110,84)(0,0)
    \put(-37,45){$\overline X~~=$}
    \put(0,0){\scalebox{.34}{\includegraphics{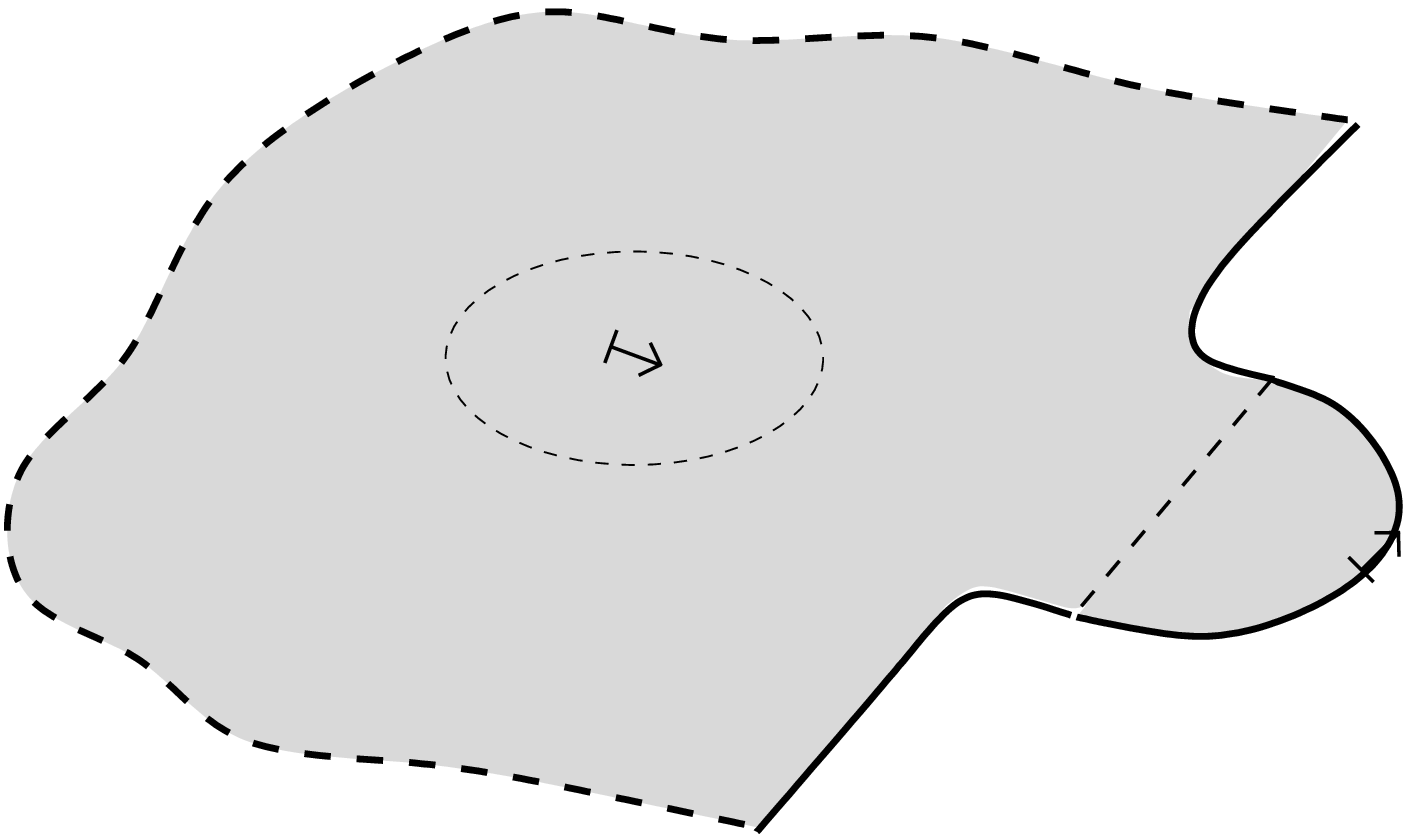}}}
    \put(125,57){\sse phys.}
    \put(120,50){\sse bnd.}
\end{picture}$$
Second, the double $\Xh$ is defined to be the orientation bundle 
${\rm Or}(\overline\X)$ over $\overline\X$,
divided by an equivalence relation that identifies the two points of 
fibres over the boundary of $\overline\X$,
$$ \begin{picture}(100,104)(0,0)
    \put(-38,59){$\Xh~~=$}
    \put(0,0){
       \put(0,0){\scalebox{.29}{\includegraphics{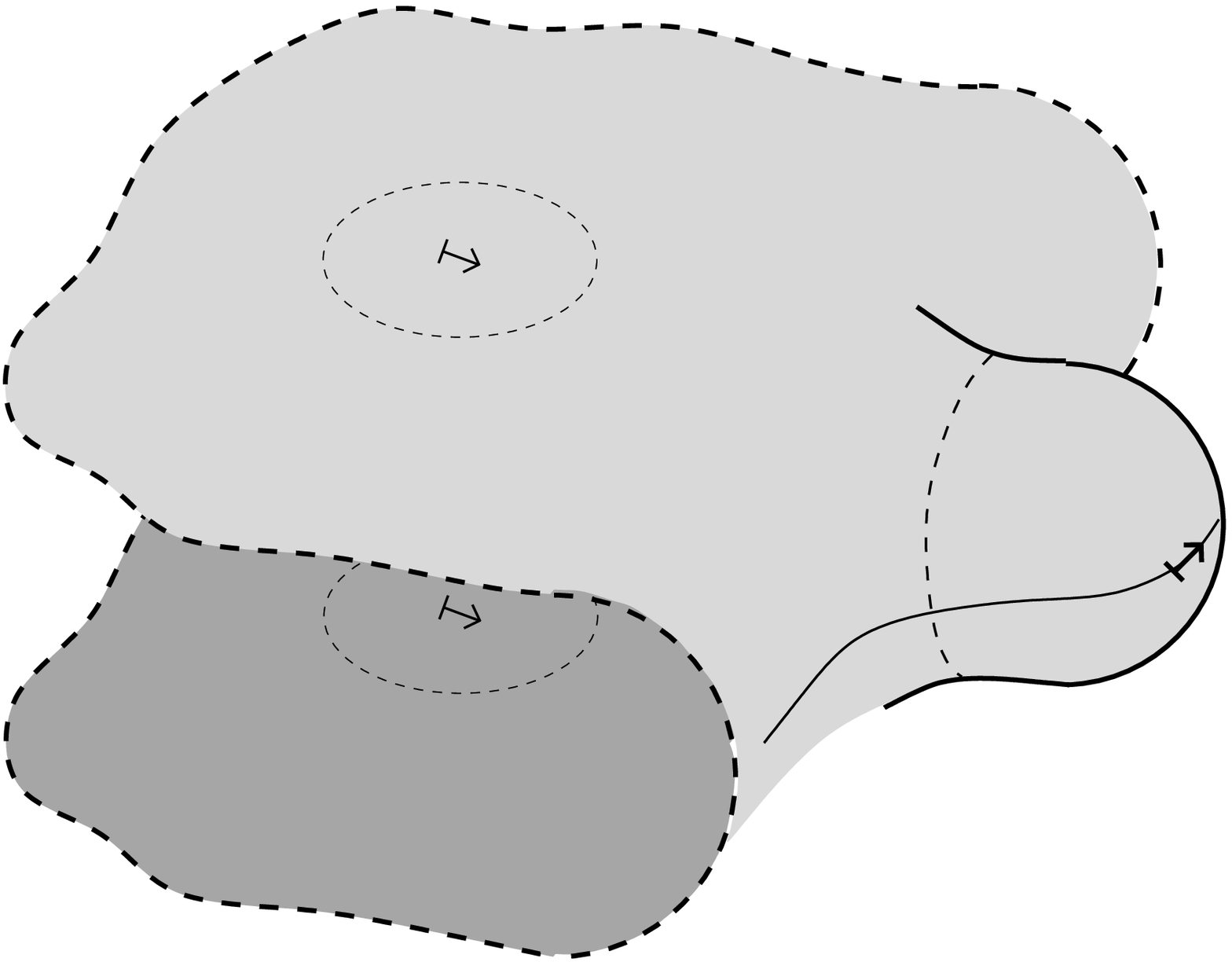}}}
       \setlength{\unitlength}{.3pt}
       \put(138,219){\sse $(B_r,-)$}
       \put(138,84){\sse $(B_l,+)$}
       \put(376,165){\sse $(A,+)$}
       \setlength{\unitlength}{1pt}
       }
\end{picture}$$
An arc embedded in the interior of $\overline\X$ has two preimages in \Xh;
we label them by $(B_l,+)$ and by $(B_r,-)$, respectively.\,\footnote{~%
  Both the world sheet \X\ (and with it $\overline X$) and the double \Xh\
  are oriented. The labelling is such that the projection from \Xh\ to
  $\overline\X$ is orientation preserving in a neighbourhood of the arc labelled
  $(B_l,+)$ and orientation reversing in a neighbourhood of the other.}
An arc on the boundary $\partial\overline X$ has one preimage in \Xh, which we
label by $(A,+)$. If a state boundary we start from is outgoing rather than 
incoming, then the corresponding plus sign in this prescription is replaced by 
a minus sign, and vice versa. There is also a natural choice
for the Lagrangian subspace,
       see appendix A.1 of \cite{fjfrs}.
In this way, \Xh\ becomes an extended surface in $\cob_\calc$.
We define the vector space $V(\X)$
as the image of \Xh\ under the functor \tftc,
\be
  V(\X) := \tftc( \Xh ) \,\in \obj(\Vectkoerper) \,.
\ee
{}From an isomorphism $\varphi\colon~\X\,{\overset\simm\longrightarrow}\,{\rm Y}$
we can construct a cobordism $\M_\varphi\colon~ \X \To{\rm Y}$ by taking the two
cylinders $\X\,{\times}\,[-1,0]$ and ${\rm Y}\,{\times}\,[0,1]$ and identifying
$\X \,{\times}\, \{ 0 \}$ with ${\rm Y} \,{\times}\, \{ 0 \}$ using $\varphi$.
Applying the \dTFT\ to this cobordism $\M_\varphi$ results in the isomorphism
\be
  \varphi^\sharp := \tftc(\M_\varphi ) :\quad V(\X)
  \overset\simm\longrightarrow V({\rm Y}) \,.
\ee
More details can again be found in \cite{TUra}.

\subsection{Cobordism for ${\rm tr}_{\rm last}$}
\label{sec:trlast}

The final ingredient in the statement of Problem \ref{prob:lem}
is the map ${\rm tr}_{\rm last}\colon~ V({\rm cut}_f(\X)) \To V(\X)$.
A precise notation would be ${\rm tr}_{\rm last}(\X,f)$ so as to keep track
of the world sheet and the parametrised cut locus, but we will use the short 
hand ${\rm tr}_{\rm last}$. It is the analogue of both ${\rm tr}_{\rm last,cl}$ 
and ${\rm tr}_{\rm last,op}$ from section \ref{sec:open-closed}, depending on 
whether $f\colon~ \rmA_\eps \To \X$ or $f\colon~ \rmS_\eps \To \X$.

By definition,
$V(\X)\eq\tftc(\Xh)$. The map ${\rm tr}_{\rm last}$ will be expressed as
\be
  {\rm tr}_{\rm last} = \tftc(\M_f) \,,~\quad {\rm where} \quad
  \M_f :~~~ \widehat{ {\rm cut}_f(\X) } \To \Xh
  \label{eq:trlast}
\ee
is a suitable cobordism. Let us start with the case that $f$ is an orientation
preserving embedding $\rmS_\eps \To \X$. Then locally the world sheet and
the cut world sheet look as follows:
\be
  \begin{picture}(235,91)(0,0)
  \put(-3,53){$\X~= $}
  \put(125,53){${\rm cut}_f(\X)~= $}
  \put(20,5){\scalebox{.5}{\includegraphics{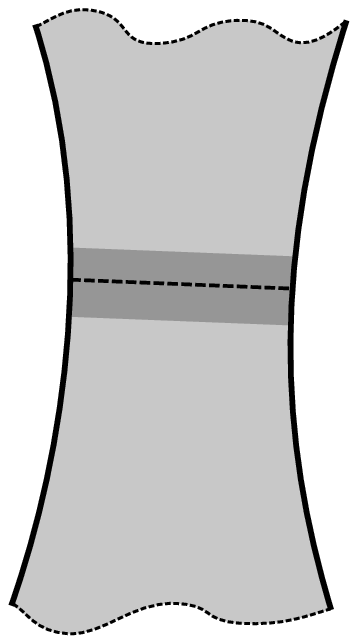}}}
  \put(180,0){
     \put(0,0){\scalebox{.5}{\includegraphics{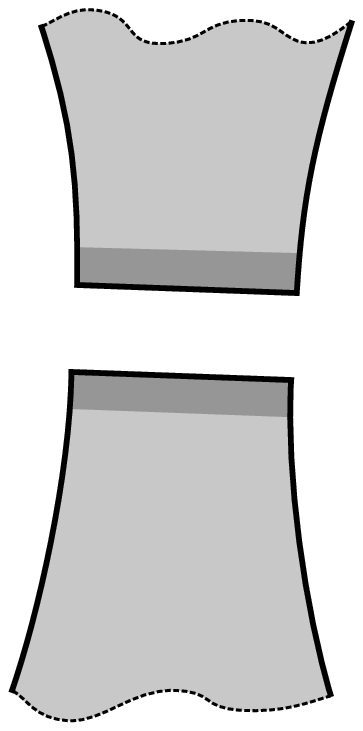}}}
     \put(45,62){\sse o-in}
     \put(44,48){\sse o-out}
     }
  \end{picture}
  \label{eq:interval-cut}
\ee
Their doubles are constructed as described in section \ref{sec:X-V(X)},
leading to
$$
  \begin{picture}(230,87)(0,0)
  \put(-5,42){$\Xh~= $}
  \put(125,42){$\widehat{ {\rm cut}_f(\X) }~= $}
  \put(20,0){\scalebox{.5}{\includegraphics{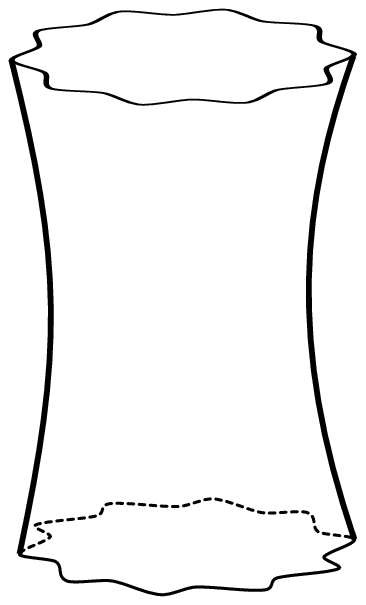}}}
  \put(180,0){
     \put(0,0){\scalebox{.5}{\includegraphics{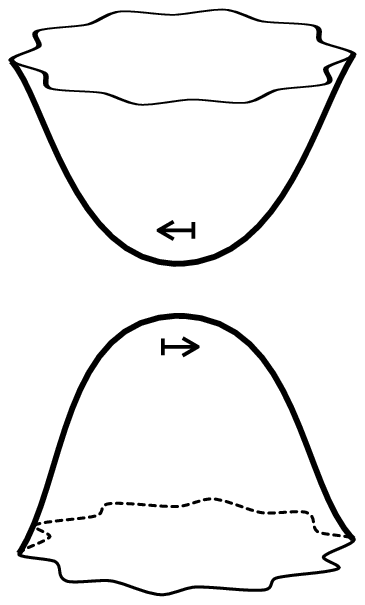}}}
     \put(14,57){\sse $(A,+)$}
     \put(14,28){\sse $(A,-)$}
     }
  \end{picture}
$$
Finally, the cobordism $\M_f$ is obtained by taking the cylinder over
$\widehat{ {\rm cut}_f(\X) }$ and identifying a disc (in fact, the disc
resulting from taking the double of the half-discs glued to the open state
boundaries) around the marked arcs at one of its ends:
\be
  \begin{picture}(250,99)(0,0)
  \put(0,48){$\M_f ~=~ \widehat{ {\rm cut}_f(\X) } \times [0,1] / {\sim}
  ~~~=$}
  \put(136,0){
     \put(0,0){\scalebox{.22}{\includegraphics{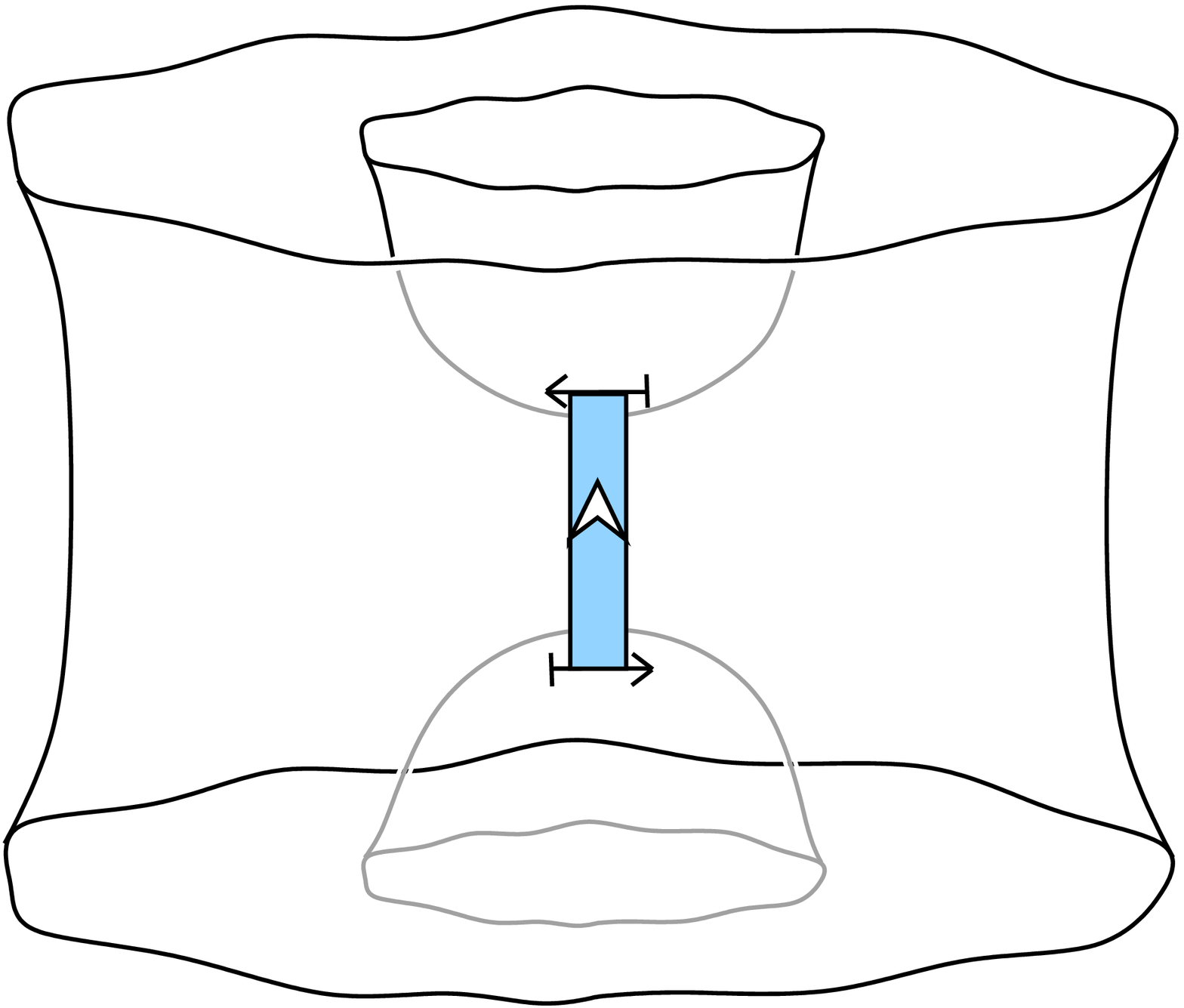}}}
     \put(48,65){\tiny $(A,+)$}
     \put(48,28){\tiny $(A,-)$}
     \put(65,45){\sse $A$}
     }
  \end{picture}
  \label{eq:trlast-open}
\ee

When cutting along a circle, we essentially have to duplicate the
above construction. Let $f\colon~ \rmA_\eps \To \X$ be an orientation
preserving embedding. Locally, the world sheet and the cut world sheet
now look as
\be
  \begin{picture}(290,105)(0,0)
  \put(-25,45){$\X~= $}
  \put(145,45){${\rm cut}_f(\X)~= $}
  \put(0,10){\scalebox{.2}{\includegraphics{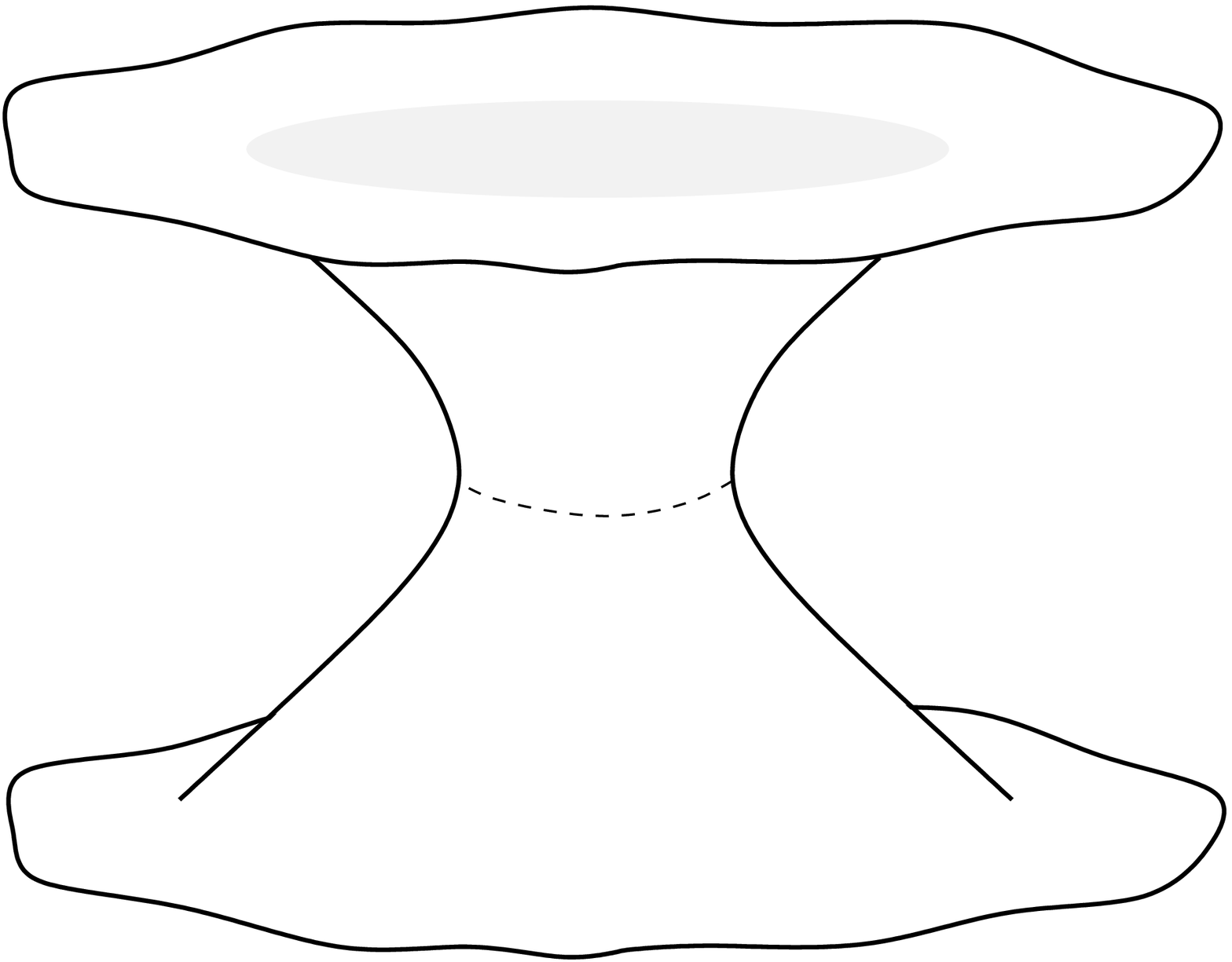}}}
  \put(200,0){
       \put(0,0){\scalebox{.2}{\includegraphics{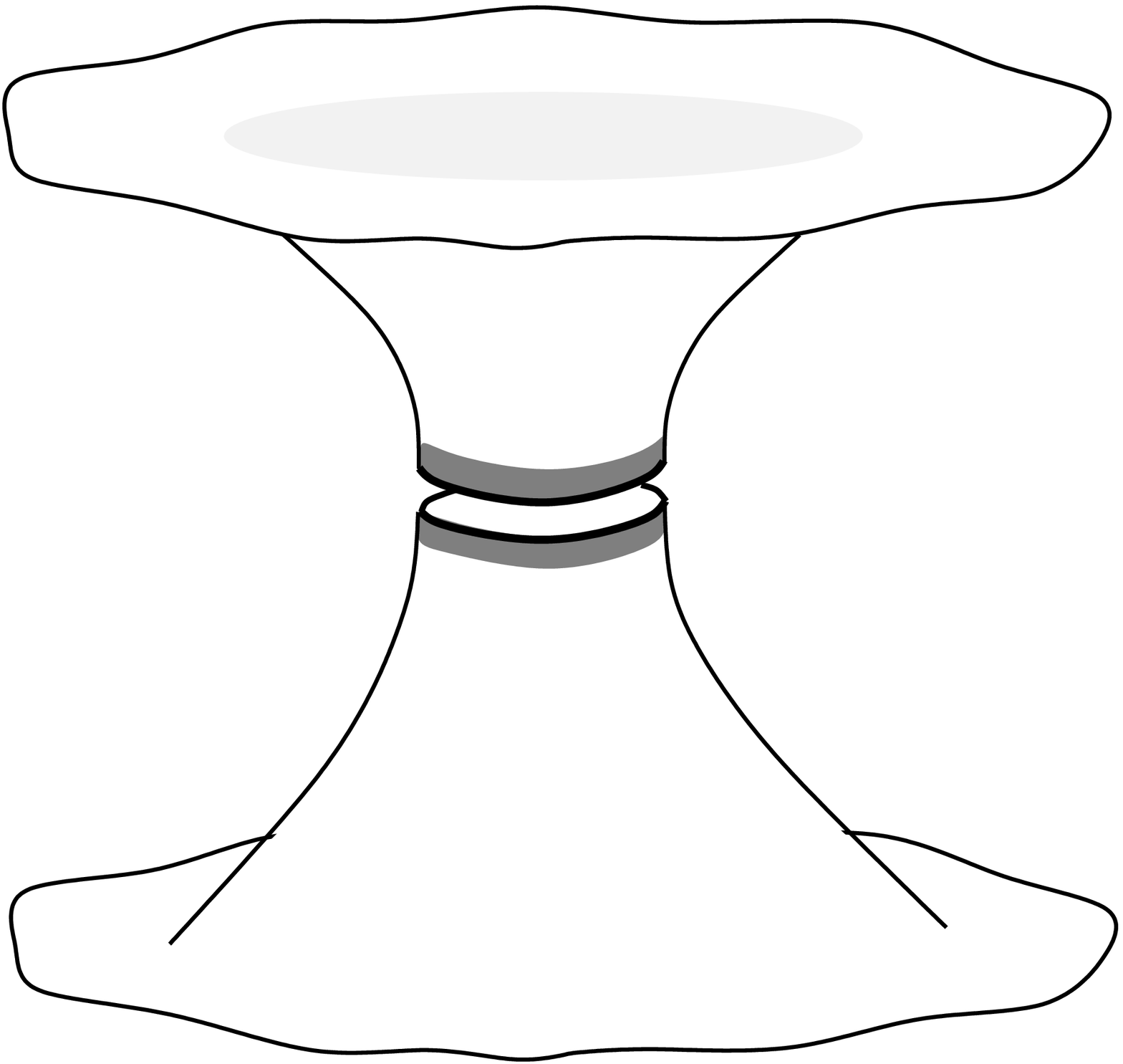}}}
       \setlength{\unitlength}{.2pt}
       \put(330,293){\sse c-in}
       \put(330,248){\sse c-out}
       \setlength{\unitlength}{1pt}
       }  
  \end{picture}
  \label{eq:cut-X-circle}
\ee
while their doubles are
\be
  \begin{picture}(290,160)(0,0)
  \put(-30,80){$\Xh~= $}
  \put(145,80){$\widehat{ {\rm cut}_f(\X) }~= $}
  \put(0,10){
       \put(0,0){\scalebox{.19}{\includegraphics{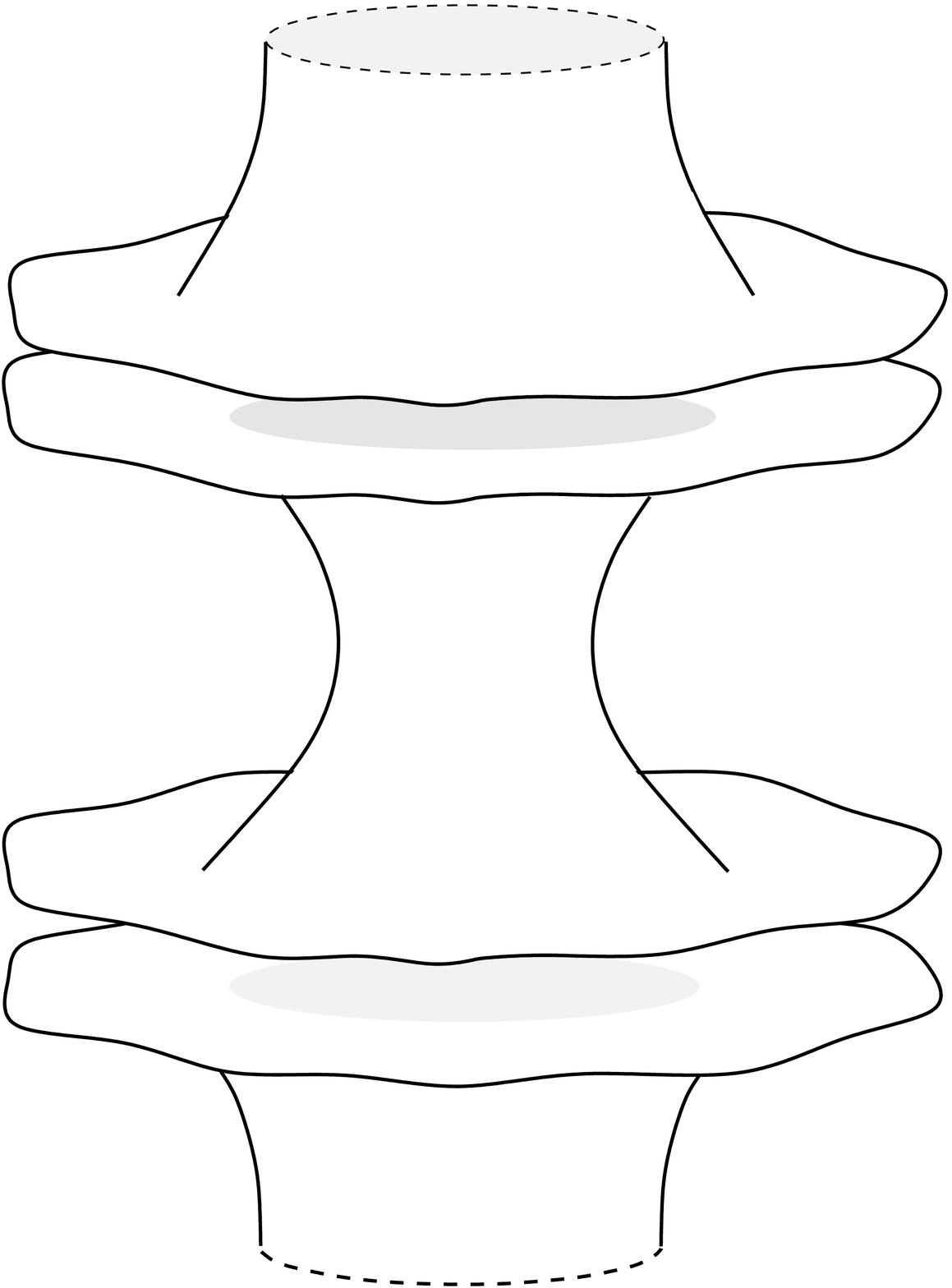}}}
       \setlength{\unitlength}{.2pt}
       \put(219,620){\sse $C$}
       \put(219,15){\sse $C'$}
       \setlength{\unitlength}{1pt}
       }
  \put(200,0){
       \put(0,0){\scalebox{.19}{\includegraphics{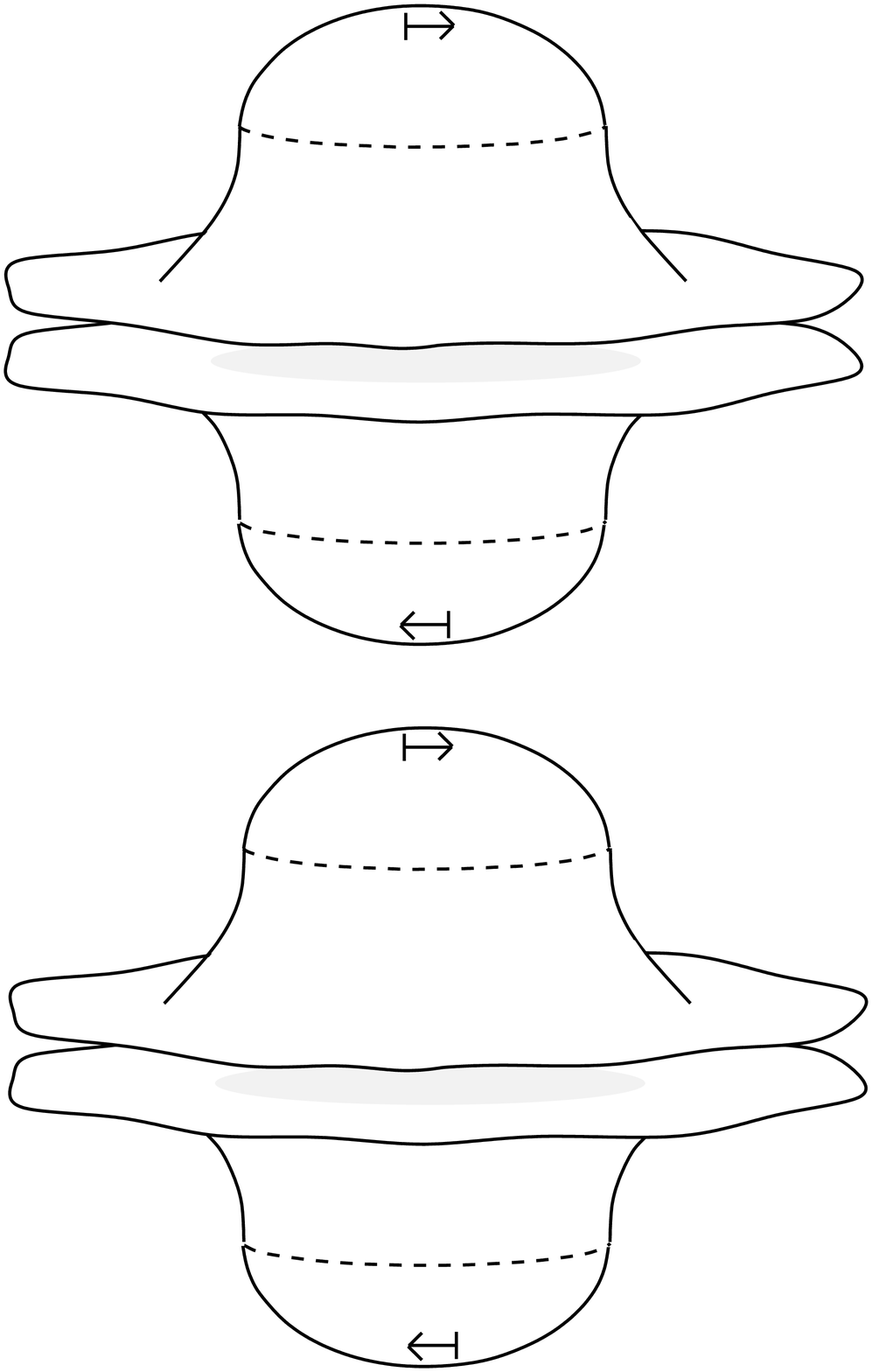}}}
       \setlength{\unitlength}{.2pt}
       \put(190,755){\tiny $(B_r,-)$}
       \put(190,459){\tiny $(B_l,+)$}
       \put(190,332){\tiny $(B_l,-)$}
       \put(190, 30){\tiny $(B_r,+)$}
       \setlength{\unitlength}{1pt}
       }
  \end{picture}
  \label{eq:cut-bulk-double}
\ee
Here for \Xh\ the circles marked $C$ and $C'$ are to be identified.
As before, the co\-bor\-dism $\M_f$ is obtained by taking the cylinder over
$\widehat{ {\rm cut}_f(\X) }$ and identifying the discs around the marked arcs
at one of its ends,
\be
  \begin{picture}(220,142)(0,0)
  \put(-10,74){$\M_f ~:=~ \widehat{ {\rm cut}_f(\X) } \times [0,1]/{\sim}~~~= $}
  \put(137,0){
       \put(0,0){\scalebox{.22}{\includegraphics{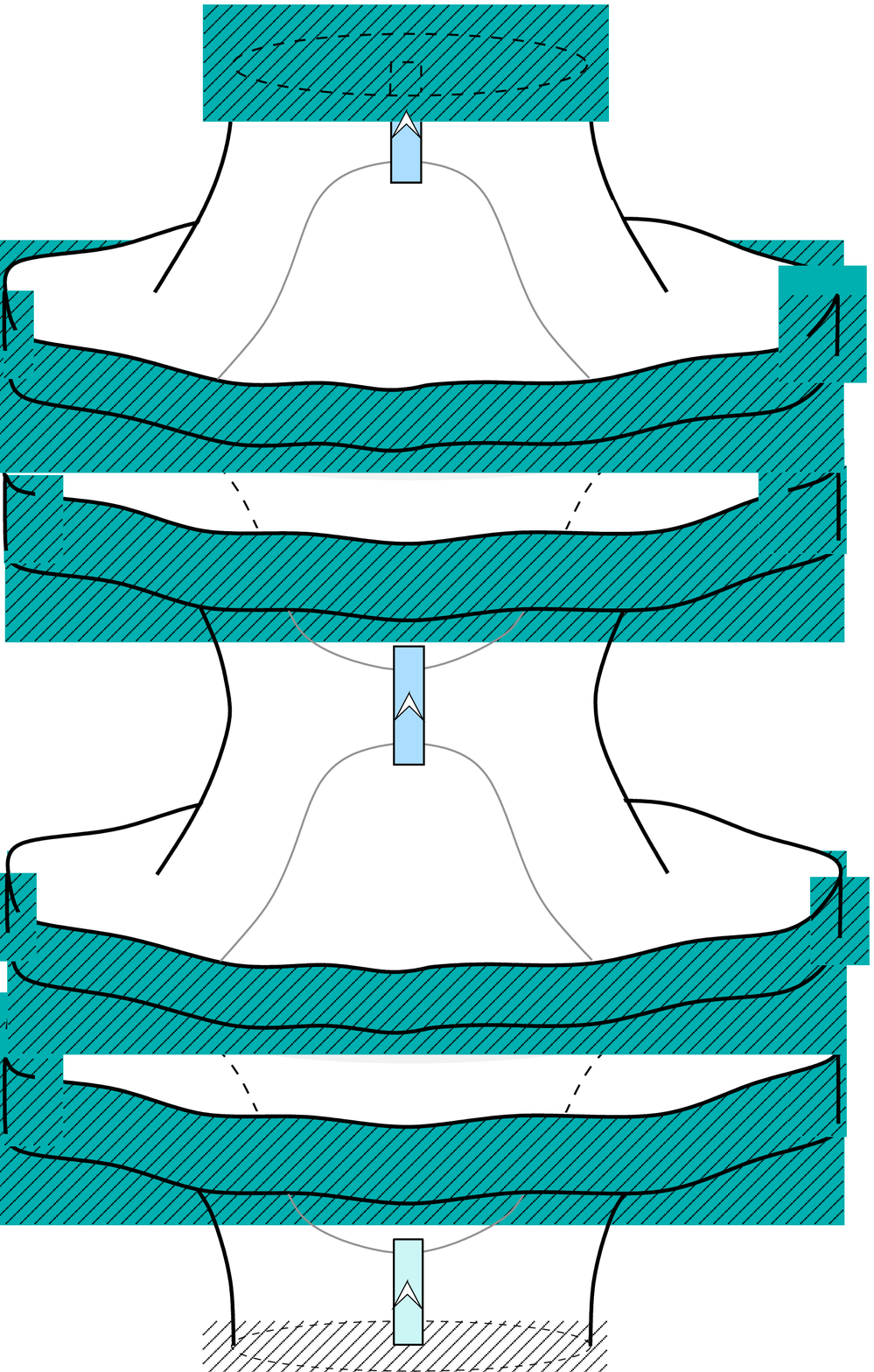}}}
       \setlength{\unitlength}{.23pt}
       \put(213,605){\tiny $B_r$}
       \put(213,320){\tiny $B_l$}
       \put(213, 33){\tiny $B_r$}
       \put( 75,631){\sse $D$}
       \put( 70,  5){\sse $D'$}
       \setlength{\unitlength}{1pt}
       }
  \end{picture}
  \label{eq:Mf-trlast-bulk}
\ee
In this picture, the discs marked $D$ and $D'$ are to be identified. Note that 
$\M_f$ is indeed a cobordism from $\widehat{ {\rm cut}_f(\X) }$ to $\Xh$.


\subsection{Frobenius algebras in $\calc$}
\label{sec:frob}

Essential for the solution to Problem \ref{prob:lem} that we will present below
is the concept of a symmetric special Frobenius algebra.  We briefly review
its definition; for more details and references consult sections 1.2 and 2.3
of \cite{C}.  A symmetric special Frobenius algebra in $\calc$ is a quintuple
$A \eq (A,m,\eta,\Delta,\eps)$, where $A \iN \obj(\calc)$ and $m,\,\eta,\,
\Delta,\,\eps$ are the morphisms of multiplication, unit, co-multiplication
and co-unit, respectively. Since $\calc$ is strict monoidal, these can be
visualised in the form of Joyal-Street type diagrams \cite{joSt5} 
(to be read from bottom to top) as follows:
$$
\quad\ m ~=
\begin{picture}(36,24)(0,17)
  \put(7,0) {\begin{picture}(0,0)(0,0)
    \scalebox{.29}{\includegraphics{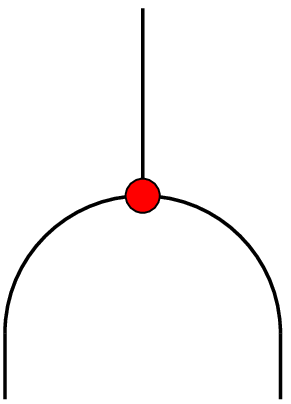}}
    \end{picture}}
  \put(3.5,-8.8)   {\sse$A$}
  \put(15.5,36.8)  {\sse$A$}
  \put(26.5,-8.8)  {\sse$A$}
\end{picture}
\quad,\qquad
\eta ~=
\begin{picture}(19,0)(0,12)
  \put(7,0) {\begin{picture}(0,0)(0,0)
    \scalebox{.29}{\includegraphics{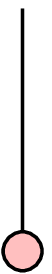}}
    \end{picture}}
  \put(5.8,25.8)   {\sse$A$}
\end{picture}
\quad,\qquad
\Delta ~=
\begin{picture}(36,0)(0,18)
  \put(7,0) {\begin{picture}(0,0)(0,0)
    \scalebox{.29}{\includegraphics{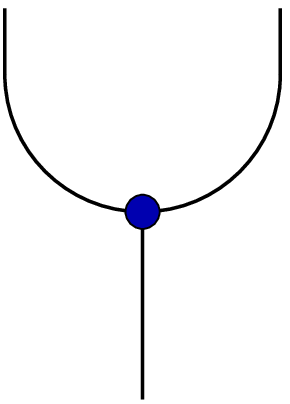}}
    \end{picture}}
  \put(3.8,36.8)   {\sse$A$}
  \put(15.2,-8.8)  {\sse$A$}
  \put(26.8,36.8)  {\sse$A$}
  \end{picture}
\quad,\qquad
\eps ~=
\begin{picture}(26,0)(0,6)
  \put(7,0) {\begin{picture}(0,0)(0,0)
    \scalebox{.29}{\includegraphics{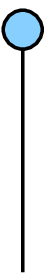}}
    \end{picture}}
  \put(5.1,-8.8)   {\sse$A$}
\end{picture}
$$
~{}\\[12pt]
In terms of this diagrammatic notation, the properties that
$m,\eta,\Delta,\eps$ must possess -- apart from (co)associativity and the 
(co)unit property --
in order that $A$ is a symmetric special Frobenius algebra look as 
$$   
\begin{picture}(253,72)
  \put(-40,58){Frobenius:}
  \put(25,0) {\begin{picture}(0,0)
    \scalebox{.29}{\includegraphics{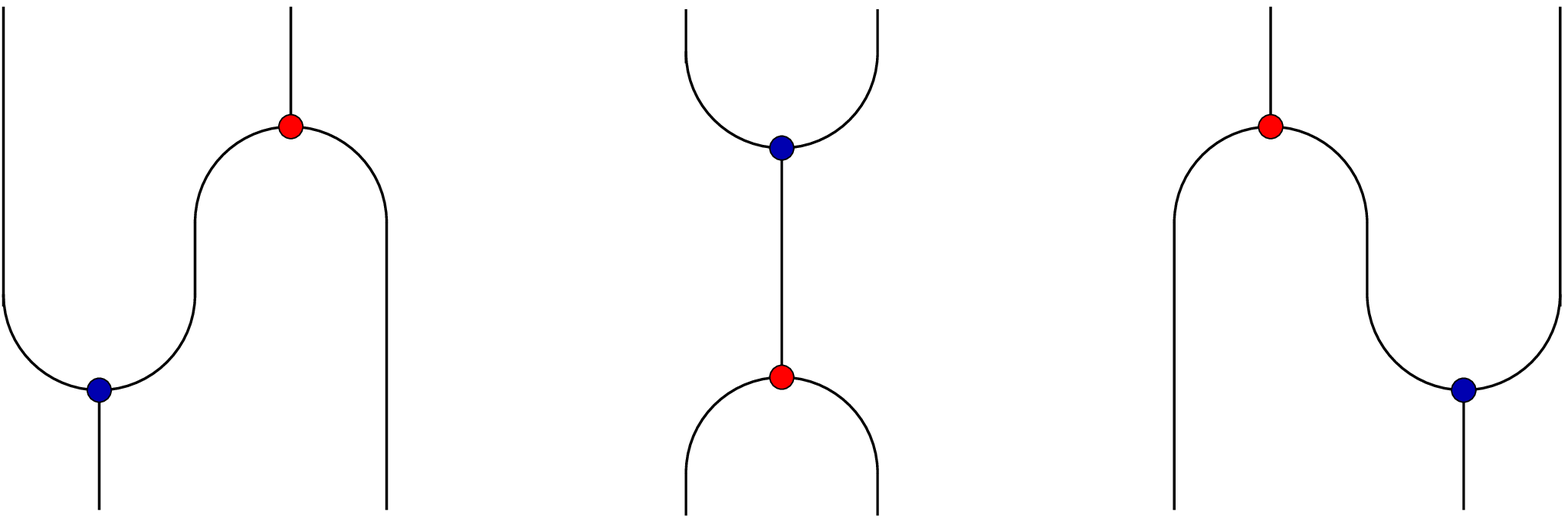}}
    \end{picture}}
  \put(3,0) {\begin{picture}(0,0)
\put(19.3,59.9)  {\sse$A$}
\put(28.7,-9.2)  {\sse$A$}
\put(51.2,59.9)  {\sse$A$}
\put(60.6,-9.2)  {\sse$A$}
\put(74.5,25)    {\small$=$}
\put(92.7,-9.2)  {\sse$A$}
\put(93.4,59.9)  {\sse$A$}
\put(113.7,-9.2) {\sse$A$}
\put(114.4,59.9) {\sse$A$}
\put(126.5,25)   {\small$=$}
\put(146.5,-9.2) {\sse$A$}
\put(157.4,59.9) {\sse$A$}
\put(178.5,-9.2) {\sse$A$}
\put(189.4,59.9) {\sse$A$}
  \end{picture}}
                 \end{picture}
$$
$$
\hspace*{2em}
\begin{picture}(155,83)(0,8)
  \put(-10,69){symmetric:}
  \put(22,0) {\begin{picture}(0,0)(0,0)
    \scalebox{.29}{\includegraphics{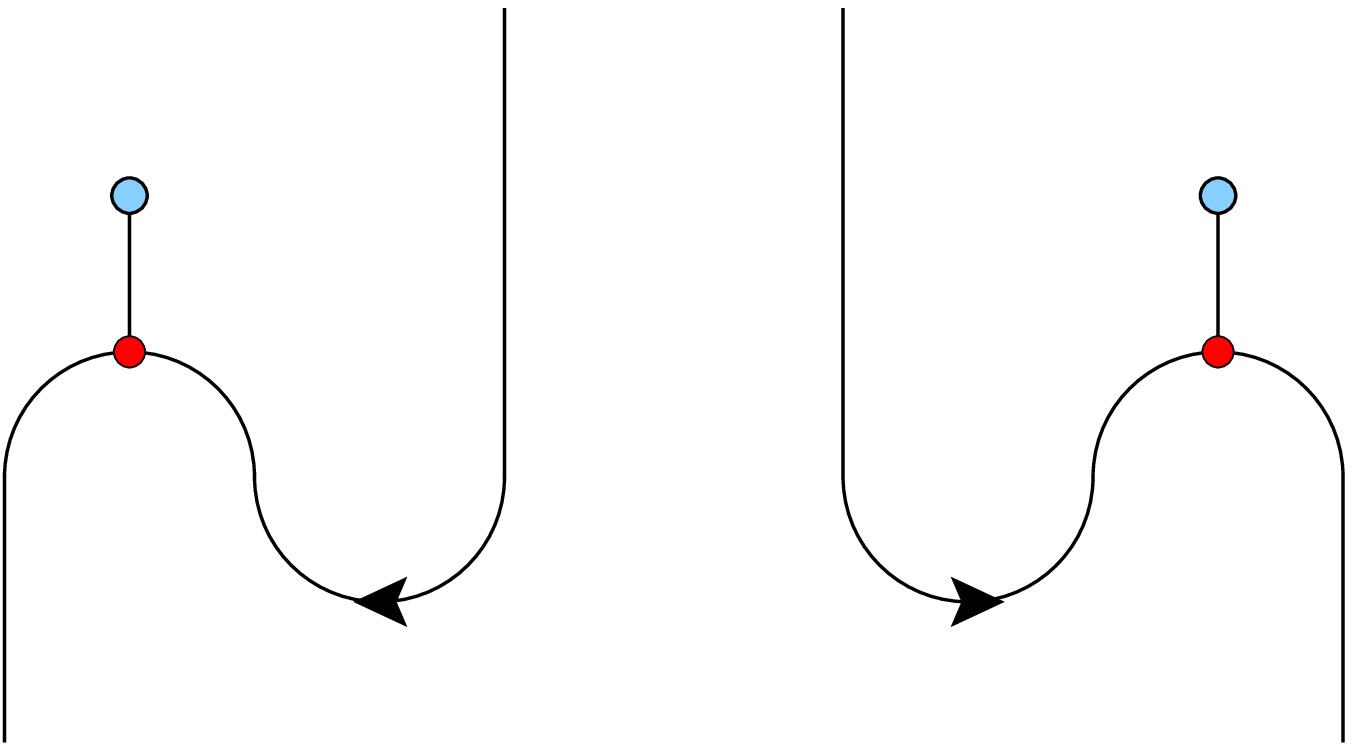}}
    \end{picture}}
\put(18.5,-9.2)  {\sse$A$}
\put(60.7,65.5)  {\sse$A^\vee$}
\put(74.9,29)    {\small$=$}
\put(89.1,65.5)  {\sse$A^\vee$}
\put(130.5,-9.2) {\sse$A$}
                 \end{picture}
\begin{picture}(99,53)
  \put(10,61){(normalised) special:}
  \put(95,20){and}
  \put(31,0) {\begin{picture}(0,0)(0,0)
    \scalebox{.29}{\includegraphics{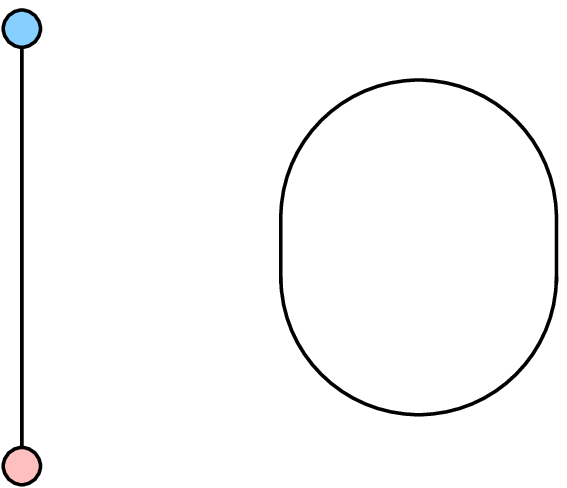}}
    \end{picture}}
\put(39.5,17)    {\small$=$}
\put(76.6,8.6)   {\sse$A$}
                 \end{picture}
\begin{picture}(99,75)
  \put(22,0) {\begin{picture}(0,0)(0,0)
    \scalebox{.29}{\includegraphics{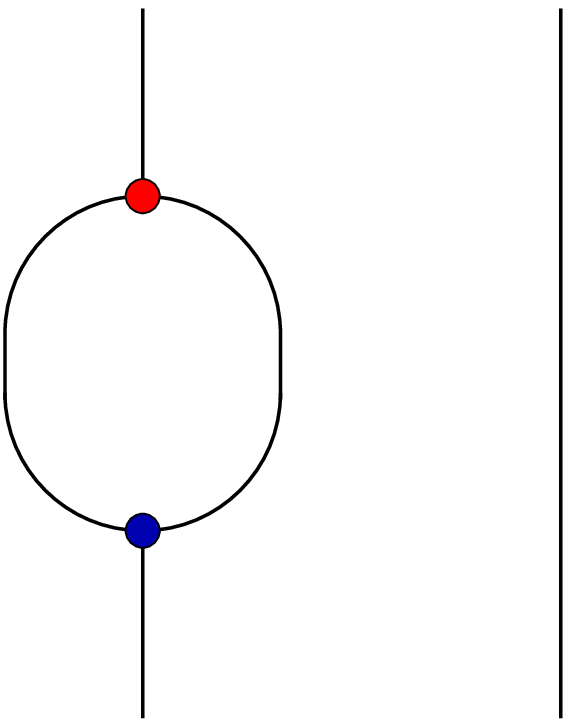}}
    \end{picture}}
\put(30.1,-9.2)  {\sse$A$}
\put(30.6,63.4)  {\sse$A$}
\put(52.5,28)    {\small$=$}
\put(64.9,-9.2)  {\sse$A$}
\put(65.4,63.4)  {\sse$A$}
                 \end{picture}
$$
~{}\\[15pt]
A Frobenius algebra is special iff $\eps \cir \eta \eq \gamma \id_\one$ and 
$m \cir \Delta \eq \gamma' \id_A$ for $\gamma,\gamma' \neq 0$. By rescaling 
$\Delta$ and $\eps$ we can always achieve $\gamma\eq\dim(A)$ and $\gamma'\eq1$. 
To emphasise that this is the choice made above we called the last property 
`normalised special'. 
Below, when we say `special', we shall always mean `normalised special'.
      
Given an algebra $A$ we can consider bimodules over $A$, and intertwiners
between them. For $A$-bimodules $B$, $B'$, we denote by $\HomAA(B,B')$ the
space of bimodule intertwiners from $B$ to $B'$.  Of particular interest
will be the bimodules $U {\otimes}^+ A {\otimes}^- V$ for
$U,V \iN \obj(\calc)$, which are defined as follows. The underlying
object is $U \oti A\oti V$, while the left and right $A$-actions are given by
$(\id_{U}\oti m\oti\id_{V}) \cir (c^{-1}_{U,A}\oti\id_A\oti\id_{V})$
and $(\id_{U}\oti m\oti\id_{V})\cir (\id_{U}\oti\id_A\oti c^{-1}_{A,V})$,
respectively (with $c_{X,Y\!}^{}$ the braiding of $\calc$).

Denote by $\II$ the label set for isomorphism classes of simple objects of
$\calc$ and by $U_i$ a representative of the class with label $i\iN \II$.
We define the object
\be
  Z(A) := \!\bigoplus_{i,j\in\II}
  \!\big( U_i{\times}\overline{U_j} \big)^{\!\oplus \tilde \Zu(A)_{ij}}
\ee
of $\calc\,{\boxtimes}\,\overline{\calc}$, with $\tilde\Zu(A)_{ij} \,{:=}\, 
\dim_\koerper^{} \HomAA(U_i{\otimes}^+\! A {\otimes}^- U_j^\vee,A)$.
The object $Z(A)$ will play the same role as the centre of
the $\complex$-algebra did in the open/closed two-dimensional
topological field theory mentioned at the end of section
\ref{sec:open-closed}.

The integers $\tilde \Zu(A)_{ij}$ and the object $Z(A)$ have a number of
interesting properties. For instance, according to Theorem 5.1 of \cite{frs1}
the matrix $\tilde \Zu(A)$ commutes with the the modular group representation
that can be constructed \cite{TUra} 
from the structural morphisms of a modular tensor category (this
property of the matrix $\tilde \Zu(A)$ was in fact first established
in the subfactor context \cite{boek1}), while as shown in
Proposition 5.3 of \cite{frs1}, the matrix $\tilde\Zu$
for the product $A\,{\atimes}\,A'$ of two symmetric special
Frobenius algebras $A$ and $A'$ is the matrix product of those of $A$ and $A'$,
$\tilde \Zu(A{\atimes}A') \eq \tilde \Zu(A)\, \tilde \Zu(A')$
(our convention for the product $\atimes$ of algebras is stated in
Proposition 3.22 of \cite{frs1}). One of the properties of the object
$Z(A)$ is given in

\begin{proposition}
\label{prop:ZA}
Let $A$ be a symmetric special Frobenius algebra in $\calc$. Then the object
$Z(A)$ of $\calc \,{\boxtimes}\, \overline{\calc}$ inherits from $A$
the structure of a \underline{commutative} symmetric Frobenius algebra.
\end{proposition}

The proof of this assertion requires several results from \cite{C};
we postpone it to appendix \ref{app:ZA-proof}.


\subsection{Construction of the correlator $\CfA(\X) \iN V(\X)$}
\label{sec:corr}

Fix a symmetric special Frobenius algebra $A$ in $\calc$.
Given a world sheet \X, we will construct an element $\CfA(\X) \iN V(\X)$
using again the \dTFT\ associated to $\calc$. As a first step
we choose a directed dual triangulation\,\footnote{~%
  That is, at every vertex precisely three edges meet, while
  faces are allowed to have an arbitrary number of edges. Also, for each edge
  a direction must be chosen, in such a way that of the three edges meeting
  at any vertex, at least one is incoming and at least one outgoing.}
$T$ of $\overline\X$ (recall from section \ref{sec:X-V(X)} that
$\overline\X$ is obtained by gluing discs and semi-discs to
the closed and open state boundaries of \X, respectively).\ We demand that
the arcs embedded in $\overline\X$ are covered by edges of $T$.

Next we introduce the {\em connecting manifold\/} $\M_A(\X,T)$. It is a
cobordism from the empty set to \Xh, defined as follows.

\smallskip\noindent
\nxt As a manifold, $\M_A(\X,T)$ is given by\,\footnote{~%
  The formula presented here applies to oriented world sheets only.
  For the unoriented case consult \cite[appendix A.1]{fjfrs}.}
\be
  \M_A(\X,T) = \overline{\X} \,{\times}\,[-1,1] /{\sim} \quad\ {\rm where}\ \
  (x,t) \,{\sim}\, (x,-t) {\rm~for~} x \iN \partial\overline{X} \,.
\ee
For example, close to a stretch of physical boundary $\M_A(\X,T)$
looks as follows:
$$
  \begin{picture}(145,82)(0,0)
  \put(0,0){\scalebox{.285}{\includegraphics{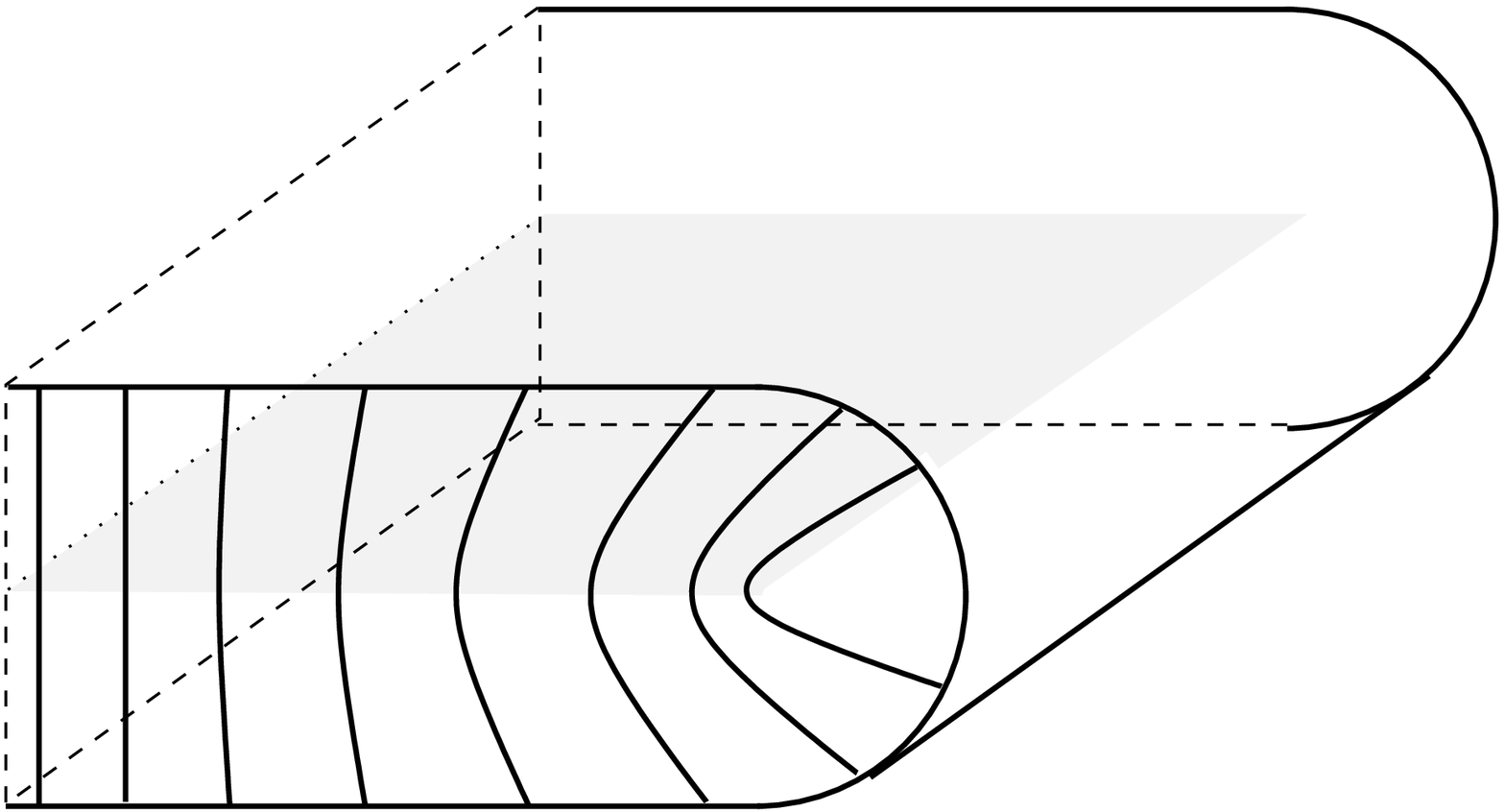}}}
  \end{picture}
$$
It is easy to see that indeed $\partial \M_A(\X,T) \eq \Xh$, as is required
if $\M_A(\X,T)$ is to yield a cobordism $\emptyset \To \Xh$. There
is also a natural embedding $\overline\X \,{\hookrightarrow}\, \M_A(\X,T)$
sending $x \iN \overline\X$ to $(x,0) \iN \M_A(\X,T)$.

\smallskip\noindent
\nxt The ribbon graph embedded in $\M_A(\X,T)$ depends on
the triangulation $T$. Thinking of $\overline{X}$ as embedded in
$\M_A(\X,T)$, ribbons labelled by $A$ are placed on the
edges of the triangulation, and the vertices are built from
the multiplication and co-multiplication morphism of $A$. Some
examples should suffice to illustrate how the ribbon graph is constructed.
Close to the boundary of $\overline{\X}$, we have
$$
  \begin{picture}(210,85)(0,0)
  \put(0,40){$\M_A(\X,T) = $}
  \put(70,0){
       \put(0,0){\scalebox{.285}{\includegraphics{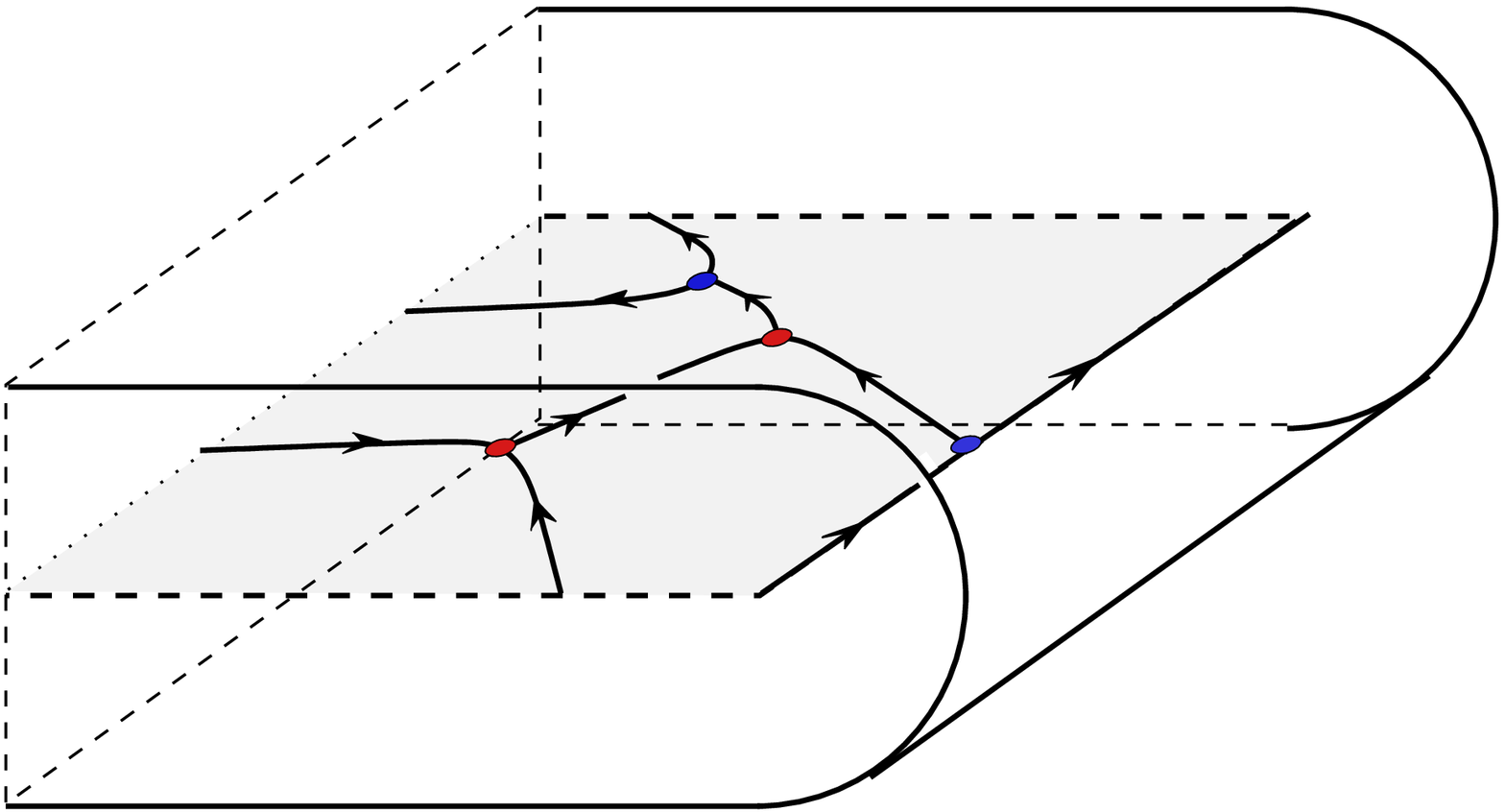}}}
       \setlength{\unitlength}{.3pt}
       \put(83,103){\tiny $A$}
       \put(369,168){\tiny $A$}
       \put(183,176){\tiny $A$}
       \setlength{\unitlength}{1pt}
       }
  \end{picture}
$$
In this picture and in the following ones we simplify the graphical
presentation by drawing ribbons as lines. Close to an arc on the boundary
(which results in this example from an incoming open state boundary of \X, 
see section \ref{sec:X-V(X)}) of $\overline{\X}$ we have
$$
  \begin{picture}(220,84)(0,0)
  \put(0,40){$\M_A(\X,T) = $}
  \put(70,0){
       \put(0,0){\scalebox{.285}{\includegraphics{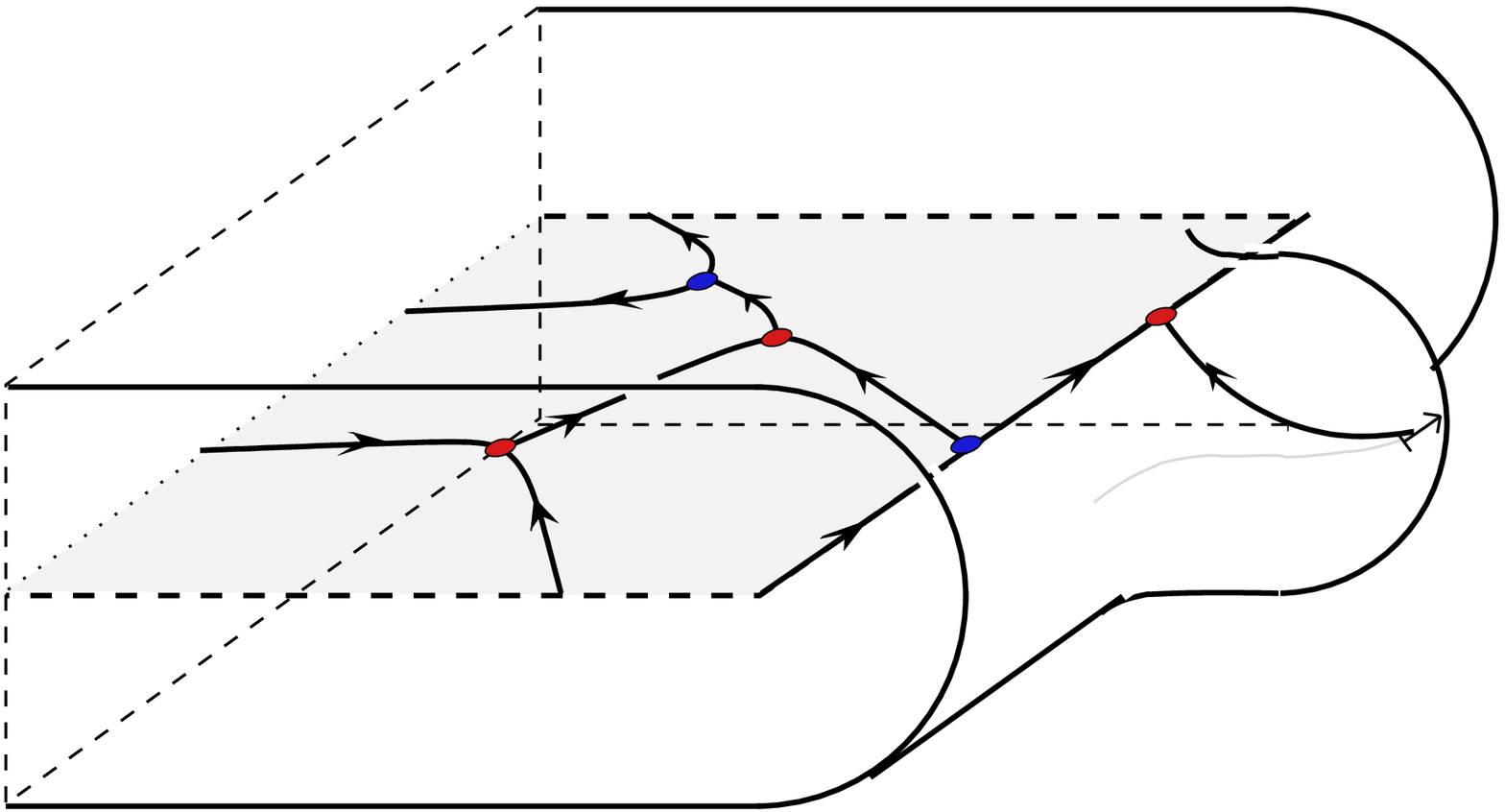}}}
       \setlength{\unitlength}{.3pt}
       \put(83,103){\tiny $A$}
       \put(369,168){\tiny $A$}
       \put(183,176){\tiny $A$}
       \put(445,133){\tiny $A$}
       \setlength{\unitlength}{1pt}
       }
  \end{picture}
$$
For an arc in the bulk, the construction is more complicated.
Consider the object
  $$
  \Omega := \bigoplus_{i \in \II} U_i \,.
  $$
Then close to an incoming (respectively outgoing) closed state boundary of $\X$
the ribbon graph in $\M_A(\X,T)$ looks as follows:
\be
  \begin{picture}(310,146)(0,0)
  \put(-20,130){a) incoming:}
  \put(150,130){b) outgoing:}
  \put(10,5){\scalebox{.333}{\includegraphics{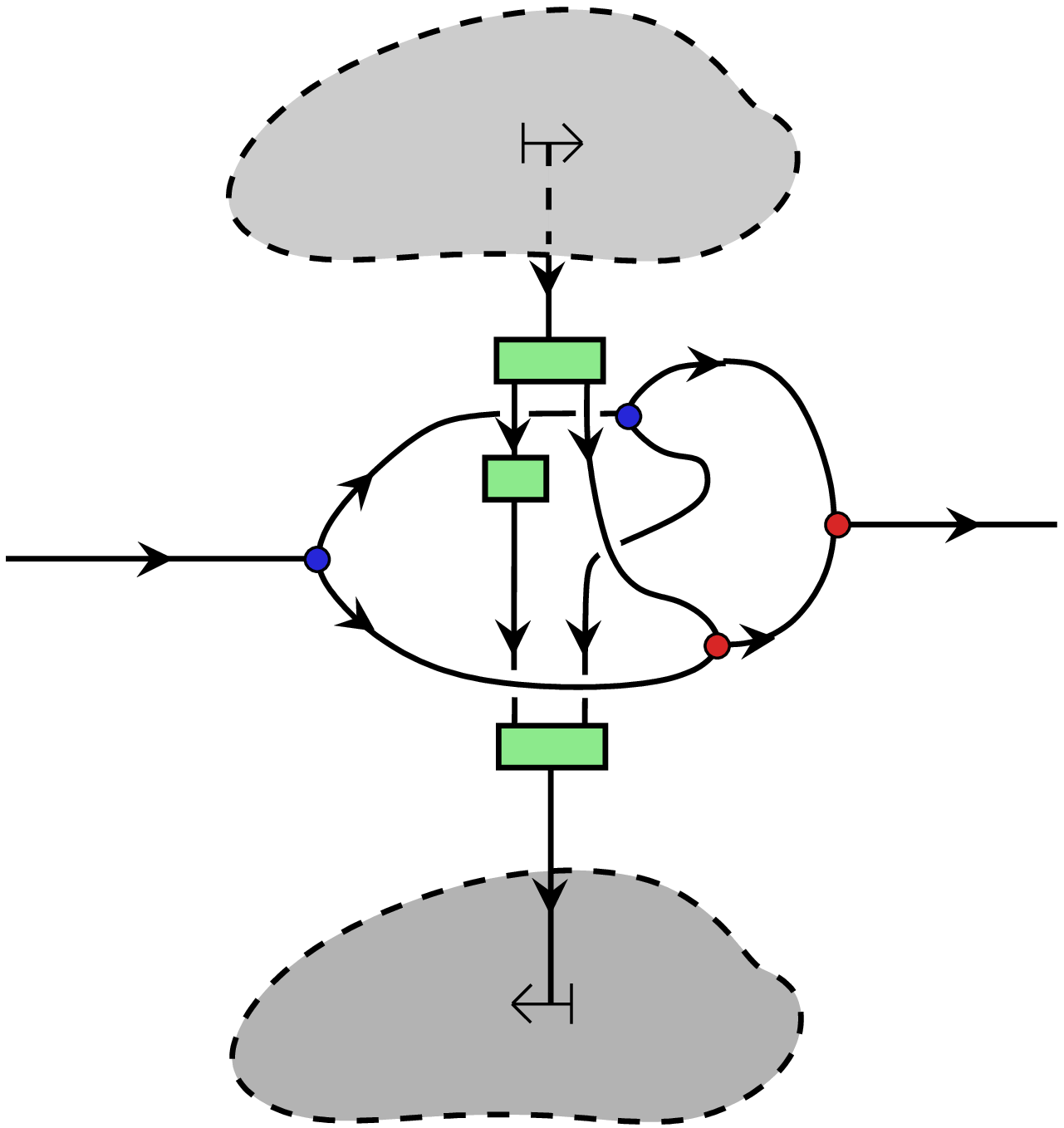}}} 
    \put(15,73){\sse$A$} \put(120,68){\sse$A$}
    \put(67.5,78.5){\tiny $\omega$} 
    \put(62,66){\sse$\Omega$}
    \put(85,68){\sse$A$}
    \put(63,100){\sse$B_l$} \put(64,38){\sse$B_r$}
    \put(59,124){\sse($B_l,+$)}
    \put(60,13){\sse($B_r,-$)}
    \put(70,92){\tiny$id$} \put(70,47){\tiny$id$}
  \put(170,5){\scalebox{.333}{\includegraphics{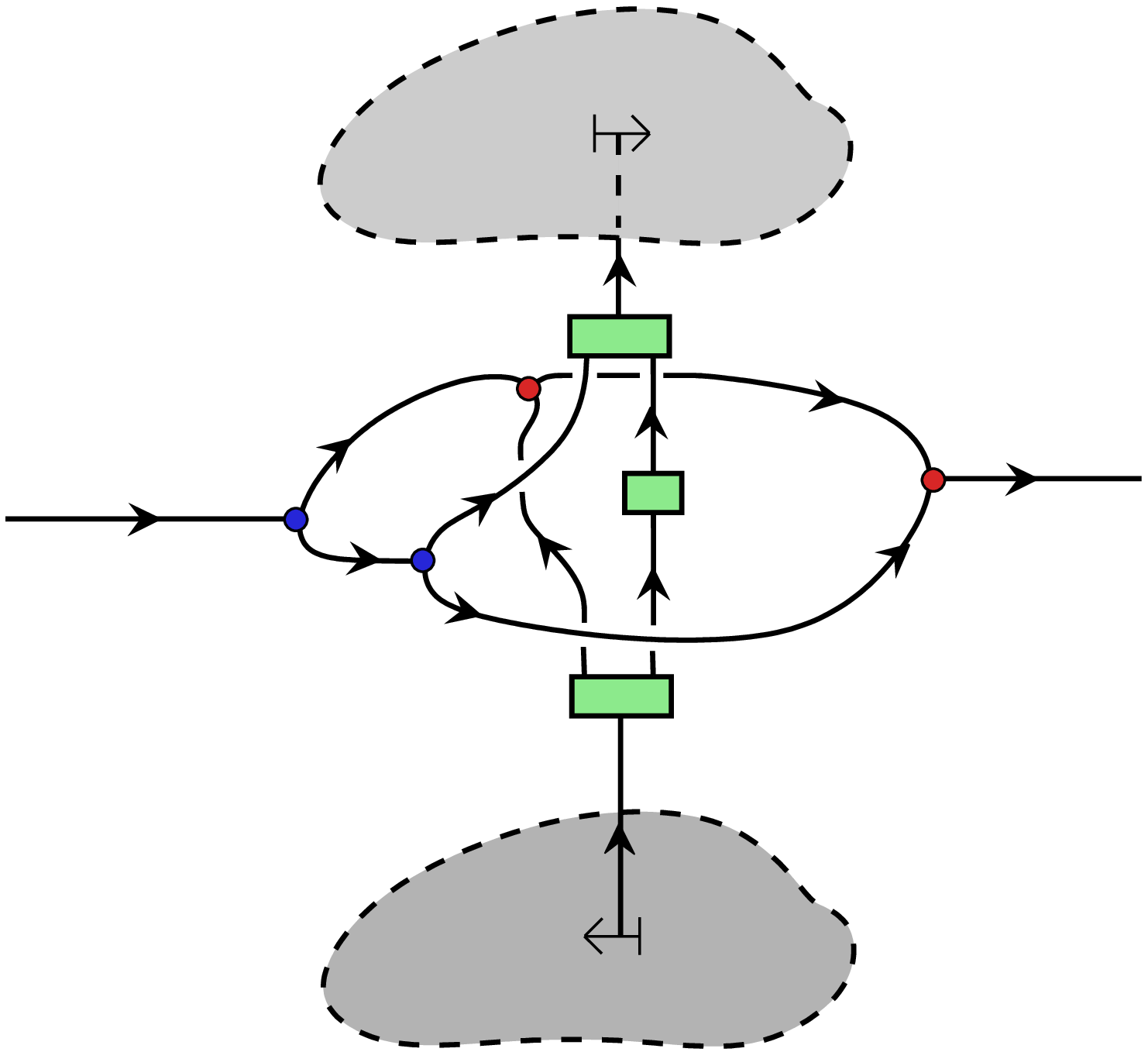}}}
    \put(175,73){\sse$A$} \put(296,69){\sse$A$}
    \put(249,73){\tiny $\omega$} 
    \put(254,60){\sse$\Omega$}
    \put(225,76){\sse$A$}
    \put(235,100){\sse$B_l$} \put(236,38){\sse$B_r$}
    \put(233,124){\sse($B_l,-$)}
    \put(233,13){\sse($B_r,+$)}
    \put(243.5,92){\tiny$id$}
    \put(243.5,47){\tiny$id$}
  \end{picture}
  \label{eq:in-out-bulk}
\ee
That is, we choose $B_l \eq B_r \eq A \oti \Omega$ (and hence it
makes sense to label the coupons joining $B_{l/r}$ and $A \oti \Omega$
by the identity morphism). The morphism $\omega \iN \End(\Omega)$ 
appearing in (\ref{eq:in-out-bulk}) is defined as
\be
  \omega := \sum_{i \in \II} \big( S_{0,0} \dim(U_i) \big)^{1/2} P_i
  \quad\ {\rm with} \quad\
  S_{0,0} = \big( \sum_{i \in\II} \dim(U_i)^2 \big)^{-1/2} .
\ee
Here $P_i \iN \End(\Omega)$ is the idempotent whose image is the subobject
$U_i$ of $\Omega$, and it is understood that we make once and for
all a choice for all square roots involved.

\smallskip\noindent
\nxt
By construction, $\tftc(\M_A(\X,T))$ is a linear map from $\koerper$ to
$V(\X)$.\ Let $T'$ be another directed dual triangulation of $\overline{\X}$.\
One can prove \cite[Proposition 3.1]{fjfrs} that
\be
  \tftc(\M_A(\X,T)) = \tftc(\M_A(\X,T'))  \,.
  \label{eq:T=T'}
\ee
The basic idea of the proof is that the properties of $A$, as depicted
in section \ref{sec:frob} are similar to the two-dimensional Matveev moves
by which one can relate any two dual triangulations of a given surface.
Because of \eqref{eq:T=T'} it makes sense to set
\be
  \CfA(\X) := \tftc(\M_A(\X,T)) 1 \in V(\X) \,,
  \label{CfA}
\ee
as the right hand side independent of the choice of $T$.


\subsection{Main theorem}

We have now gathered all ingredients needed to state our main result.

\begin{theorem}
\label{thm:solve}
Let $A$ be a symmetric special Frobenius algebra in $\calc$.
Then a solution to Problem \ref{prob:lem} is obtained
by setting $H_{\rm op} \eq A$, $H_{\rm cl} \eq Z(A)$
and taking the assignment $\X \,{\mapsto}\, \CfA(\X)$ as in \eqref{CfA}.
\end{theorem}

By formulating the theorem in terms of $\CfA$ it is understood that
in the prescription below Problem \ref{prob:lem} one chooses
$B_l \eq B_r \eq A \oti \Omega$.

Theorem \ref{thm:solve} has been proven, in a formulation more oriented
towards using simple objects instead of the direct sum $\Omega$, in
Theorems 2.1, 2.6 and 2.9 of \cite{fjfrs}. Let us give a sketch of the proof.
First, the requirement {\bf(A1)} in Problem \ref{prob:lem} follows
form the triangulation independence of $\CfA(\X)$. Given an isomorphism
$\varphi\colon~ \X \To {\rm Y}$, we can transport a triangulation $T$ from
$\overline\X$ to a triangulation $T'$ of $\overline{\rm Y}$. It is then not
difficult to see that $\M_\varphi \cir \M_A(\X,T)$ and $\M_A({\rm Y},T')$
are in the same homotopy class of cobordisms. It follows that
\be
  \begin{array}{ll}
  \varphi^\sharp \Cf(\X) \!\!\!\!&= \tftc(\M_\varphi) \circ
  \tftc(\M_A(\X,T)) = \tftc( \M_\varphi \cir \M_A(\X,T) )\\{}\\[-.7em]
  &= \tftc( \M_A({\rm Y},T') ) = \Cf({\rm Y}) \,.
  \end{array}
\ee

For an embedding $f\colon~ \rmS_\eps \To \X$ also property {\bf(A2)} is easy to
check. Just recall that cutting a world sheet along an interval looks locally
as in \eqref{eq:interval-cut} and that the cobordism for ${\rm tr}_{\rm last}$
takes the form \eqref{eq:trlast-open}. Combining with the prescription to
construct the connecting manifold, the desired equality
$\CfA(X)\eq{\rm tr}_{\rm last}( {\rm cut}_f(\X) )$ amounts to verifying that
$$
  \begin{picture}(240,111)(0,0)
  \put(-25,55){$\tftc\Big($}
  \put(95,55){$\Big) \quad = \quad \tftc\Big($}
  \put(256,55){$\Big)$}
  \put(0,0){
       \put(0,0){\scalebox{.24}{\includegraphics{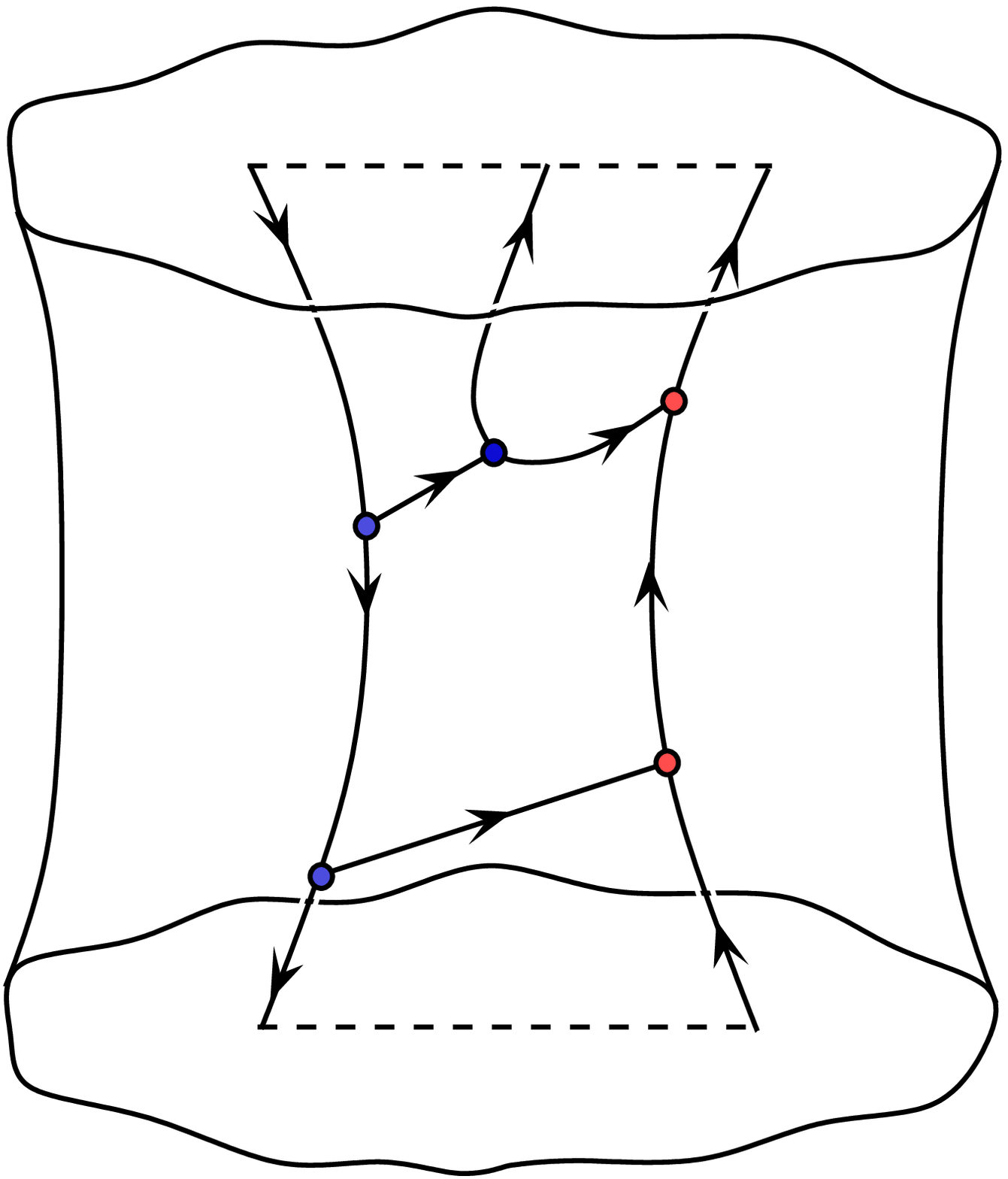}}}
       \setlength{\unitlength}{.25pt}
       \put(107,227){\tiny $A$}
       \put(178,145){\tiny $A$}
       \put(248,222){\tiny $A$}
       \put(202,341){\tiny $A$}
       \setlength{\unitlength}{1pt}
       }
  \put(160,0){
       \put(1,1){\scalebox{.24}{\includegraphics{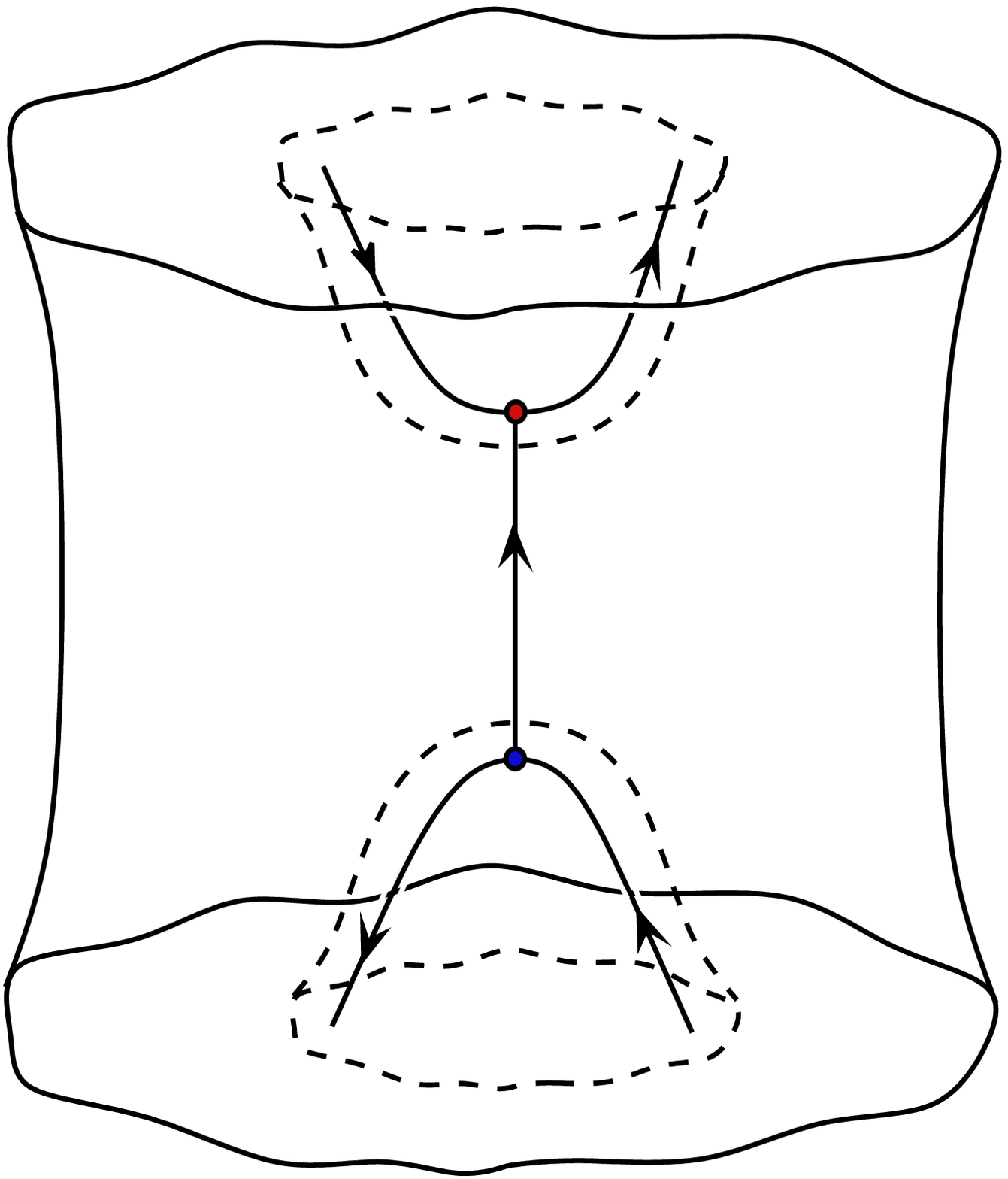}}}
       \setlength{\unitlength}{.25pt}
       \put(203,220){\sse $A$}
       \put(131,380){\sse $A$}
       \put(230, 60){\sse $A$}
       \setlength{\unitlength}{1pt}
       }
  \end{picture}
$$
That this equality indeed holds can be shown by making use of the
properties of $A$ given in section \ref{sec:frob} \cite{fjfrs}.

For an embedding $f\colon~ \rmA_\eps \To \X$ the geometry
is more involved. The proof of property {\bf(A2)} for this case
is sketched in appendix \ref{app:bulk-fact}.


\section{Comments}

We have argued in section \ref{sec:motivate}
that solving Problem \ref{prob:lem} is equivalent to constructing
an open/closed CFT whose underlying chiral symmetry is a given
rational vertex algebra. Section \ref{sec:algebra} was devoted to
the precise formulation of Problem \ref{prob:lem} and to proving
that solutions can be obtained from symmetric special Frobenius algebras
in the relevant modular tensor category $\calc$.
To set these considerations into perspective, let us mention a number of
issues that we did not address in the discussions above.
\\[3.4pt]
\nxt From a
quantum field theoretic point of view it is natural to allow
different physical boundary segments to carry different boundary
conditions. In the algebraic setting, boundary conditions are
left modules of the algebra $A$, see \cite{frs1}.
\\[3.4pt]
\nxt We have {\em not\/} claimed that {\em every\/} solution to Problem
\ref{prob:lem} gives rise to a symmetric special Frobenius algebra in
$\calc$. To address this point is an obvious next step in our future work.
\\[3.4pt]
\nxt Another question is whether different symmetric special Frobenius algebras
in $\calc$ yield different CFTs.\ One finds that an oriented open/closed CFT
together with {\em all\/} its symmetry preserving boundary conditions is
encoded not by a single such algebra, but by a {\em Morita class\/} of
symmetric special Frobenius algebras, or, equivalently, by an appropriate
{\em module category\/} over $\calc$, see section 4.1 of \cite{frs1} for 
details and references. On the other hand, an oriented open/closed CFT with 
{\em one\/} preferred boundary condition -- the situation discussed in the 
present article -- corresponds to an {\em isomorphism\/} class of symmetric 
(normalised) special Frobenius algebras.
\\[3.4pt]
\nxt We have only treated open/closed CFTs on oriented
world sheets. Results similar to those discussed here can
also be obtained for {\em unoriented\/} world sheets, in which case
suitable equivalence classes of so-called Jandl algebras are the relevant
algebraic structure, see \cite{frs2} and section 11 of \cite{us}. 
\\[3.4pt]
\nxt In another approach \cite{bppz2,pezu5} to the algebraic
part of the construction of a CFT it is argued that the structure
of a weak Hopf algebra can be extracted from a rational 
open/closed CFT. The connection to the present approach is via
module categories, see section 4 of \cite{ffrs4} for more details.

\smallskip

As a closing remark, let us emphasise the following point.
We have seen that algebras in braided monoidal categories appear naturally in
the study of CFTs. Various physical quantities on the CFT side correspond to
a standard construction on the algebraic side. Boundary conditions versus
modules over an algebra provide just one example, but there are many more; 
see e.g.\ the dictionary in section 7 of \cite{us}. One can thus try and use 
CFT as a guiding principle in generalising other aspects of algebra
from vector spaces to a truly braided setting.

The centre of an algebra $A$ provides a nice illustration of this idea. In the
braided setting one must consider two centres, the left centre $C_l$ and the
right centre $C_r$ \cite{vazh,ostr,C}. Both are subobjects of $A$, and they
coincide if the braiding is symmetric.  Comparing the example at the end
of section \ref{sec:open-closed} and the statement of Theorem \ref{thm:solve},
another natural object to consider is $Z(A)$, which 
unlike $A$ and $C_{l/r}$ is not an object of $\calc$,
but of $\calc \,{\boxtimes}\, \overline{\calc}$.\ Actually
this subsumes the notion of left and right centre in the sense that both
$C_l \,{\times}\, \one$ and $\one \,{\times}\, C_r$ are subobjects of $Z(A)$.
Further, there is a notion of a centre $\mathcal{Z}$ of a monoidal category, 
which for a modular tensor category $\calc$ is given by
$\mathcal{Z}(\calc) \,{\cong}\, \calc \,{\boxtimes}\, \overline{\calc}$
\cite{muge9}. Thus for a symmetric special Frobenius algebra
$A$ in a modular tensor category $\calc$, the object $Z(A)$
is a commutative symmetric Frobenius algebra in the centre $\mathcal{Z}(\calc)$.

Other examples of structures and results arising when generalising algebra
to the braided setting, and which become tautological when the
braiding is symmetric, can be found in \cite{C}. In fact, the issues
investigated in \cite{C} were motivated from problems in CFT, too.


\appendix

\section{Proof of Proposition \ref{prop:ZA}}\label{app:ZA-proof}

In the following we freely use notation from \cite{C}.
We will show that, as objects in $\calc \,{\boxtimes}\, \overline{\calc}$,
\be
    Z(A) \,\cong\, C_l( (A{\times}\one) \oti T_\calc ) \,,
    \label{eq:ZA-Cl-iso}
\ee
where $T_\calc\,{\cong}\,\bigoplus_i U_i \,{\times}\, \overline{U_i}$ as an
object of $\calc \,{\boxtimes}\, \overline{\calc}$ (see sections 2.4 and 6.3
of \cite{C}) and $C_l(\cdot)$ denotes the left centre.
Since both $A$ and $T_\calc$ are symmetric special Frobenius
algebras, so is $(A{\times}\one) \oti T_\calc$. Then, by Proposition 2.37(i) of
\cite{C}, the left centre $C_l( (A{\times}\one) \oti T_\calc )$ is a
commutative symmetric Frobenius algebra. Together with
the isomorphism (\ref{eq:ZA-Cl-iso}) this proves Proposition \ref{prop:ZA}.

To establish the isomorphism (\ref{eq:ZA-Cl-iso})
we combine Propositions 3.14(i) and 3.6 of \cite{C} so as to find
\be
  \begin{array}{l}
  \Hom(C_l((A\,{\times}\,\one)\oti T_\calc) , U_p \,{\times}\, \overline{U_q} )
  \\{}\\[-.5em]
  \qquad\qquad\qquad\quad \cong
  \Hom(E^l_{A\,{\times}\,\one}( T_\calc ) , U_p \,{\times}\, \overline{U_q} )
  \\{}\\[-.5em]
  \qquad\qquad\qquad\quad \cong
     \HomAA((A\,{\times}\,\one) {\otimes}^- T_\calc ,
  (A\,{\times}\,\one) {\otimes}^+ (U_p \,{\times}\, \overline{U_q}) )
  \\{}\\[-.5em]
  \qquad\qquad\qquad\quad \cong
  \bigoplus_{k\in\II}
  \HomAA(A {\otimes}^- U_k, A {\otimes}^+ U_p) \otimes_\koerper
  \Hom(\overline{U_k},\overline{U_q}) \,.
  \end{array}
\ee
Together with \cite[Lemma 2.2]{frs4} and \cite[Theorem 5.23(iii)]{frs1}
this results in
\be
  \begin{array}{l}
  \dim\Hom(C_l((A\,{\times}\,\one)\oti T_\calc ), U_p\,{\times}\,\overline{U_q})
  \\{}\\[-.5em]
  \qquad\qquad\qquad\quad = \dim\HomAA(A {\otimes}^- U_q, A {\otimes}^+ U_p)
  \qquad\qquad\qquad\qquad
  \\{}\\[-.5em]
  \qquad\qquad\qquad\quad = \dim\HomAA(U_p^\vee {\otimes}^+ A {\otimes}^- U_q,A)
  \\{}\\[-.5em]
  \qquad\qquad\qquad\quad = \tilde \Zu(A)_{\bar p, \bar q} = \tilde \Zu(A)_{p,q}
  \\{}\\[-.5em]
  \qquad\qquad\qquad\quad = \dim \Hom(Z(A) , U_p\,{\times}\,\overline{U_q}) \,.
  \end{array}
\ee
Since $U_p \,{\times}\, \overline{U_q}$ are representatives for the isomorphism
classes of simple objects of the semisimple category $\calc\,{\boxtimes}\,
\overline{\calc}$, this shows that $Z(A)$ and $C_l((A{\times}\one)\oti T_\calc)$
are indeed isomorphic as objects in $\calc \,{\boxtimes}\, \overline{\calc}$,
and thus completes the proof.


\section{Bulk factorisation}\label{app:bulk-fact}

Here we present some details of the proof that the construction
presented in section \ref{sec:corr} solves the requirement
{\bf(A2)} also when cutting along a circle, i.e.\ for an embedding
$f\colon~ \rmA_\eps \To \X$. To this end we start
from the situation depicted in (\ref{eq:cut-X-circle}).
Close to the image of the embedding $f\colon~ \rmA_\eps \To \X$
the world sheet $\X$ looks as follows.
$$
  \begin{picture}(170,117)(0,0)
  \put(0,54){$\X ~~= $}
  \put(38,0){\scalebox{.24}{\includegraphics{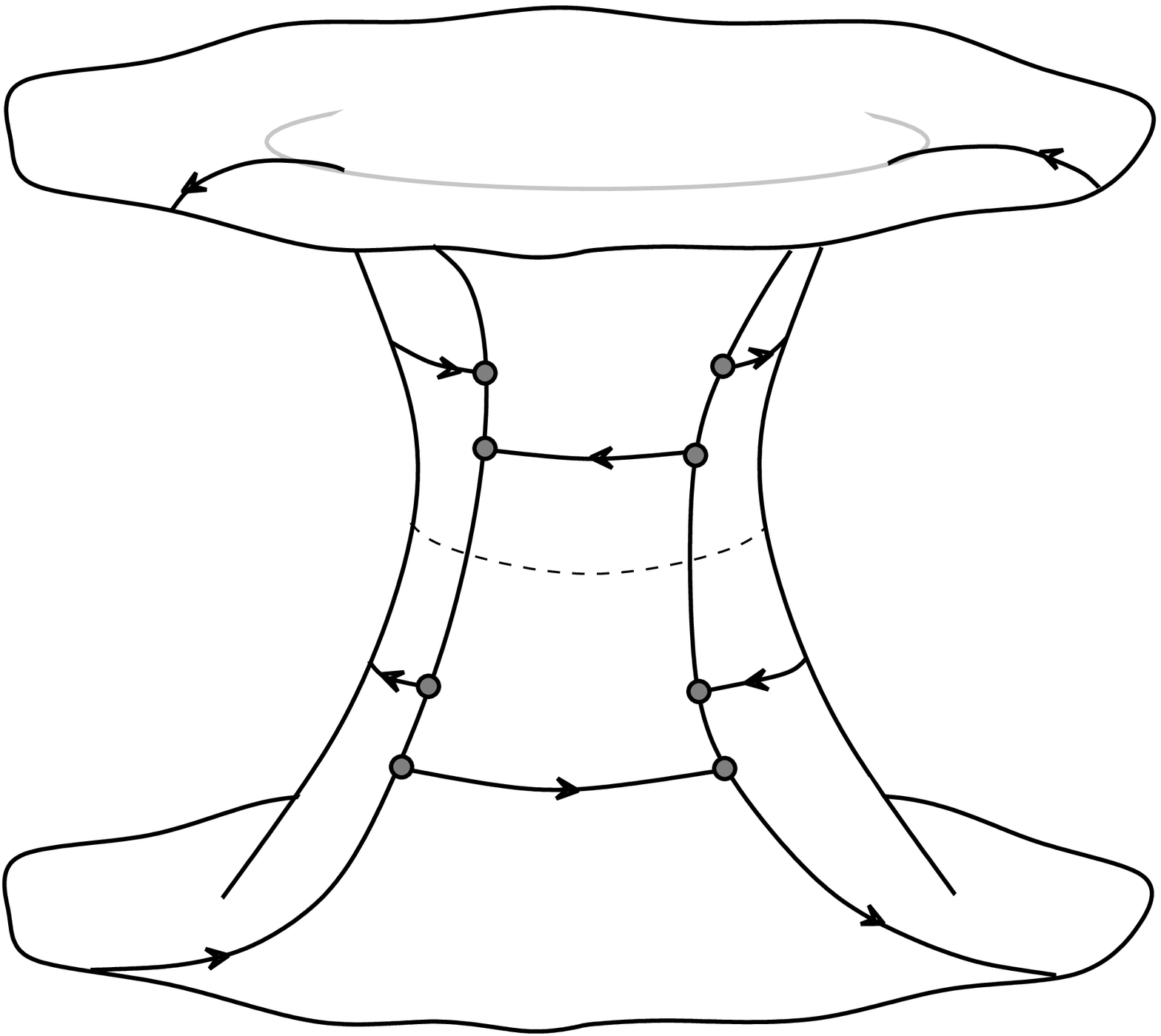}}}
  \end{picture}
$$
Here we have also indicated the dual triangulation $T$ we will use.
The corresponding fragment of the connecting manifold takes the form
\be
  \begin{picture}(160,123)(0,0)
  \put(-22,57){$\M_A(\X,T) ~~= $}
  \put(50,-1){
       \put(0,0){\scalebox{.28}{\includegraphics{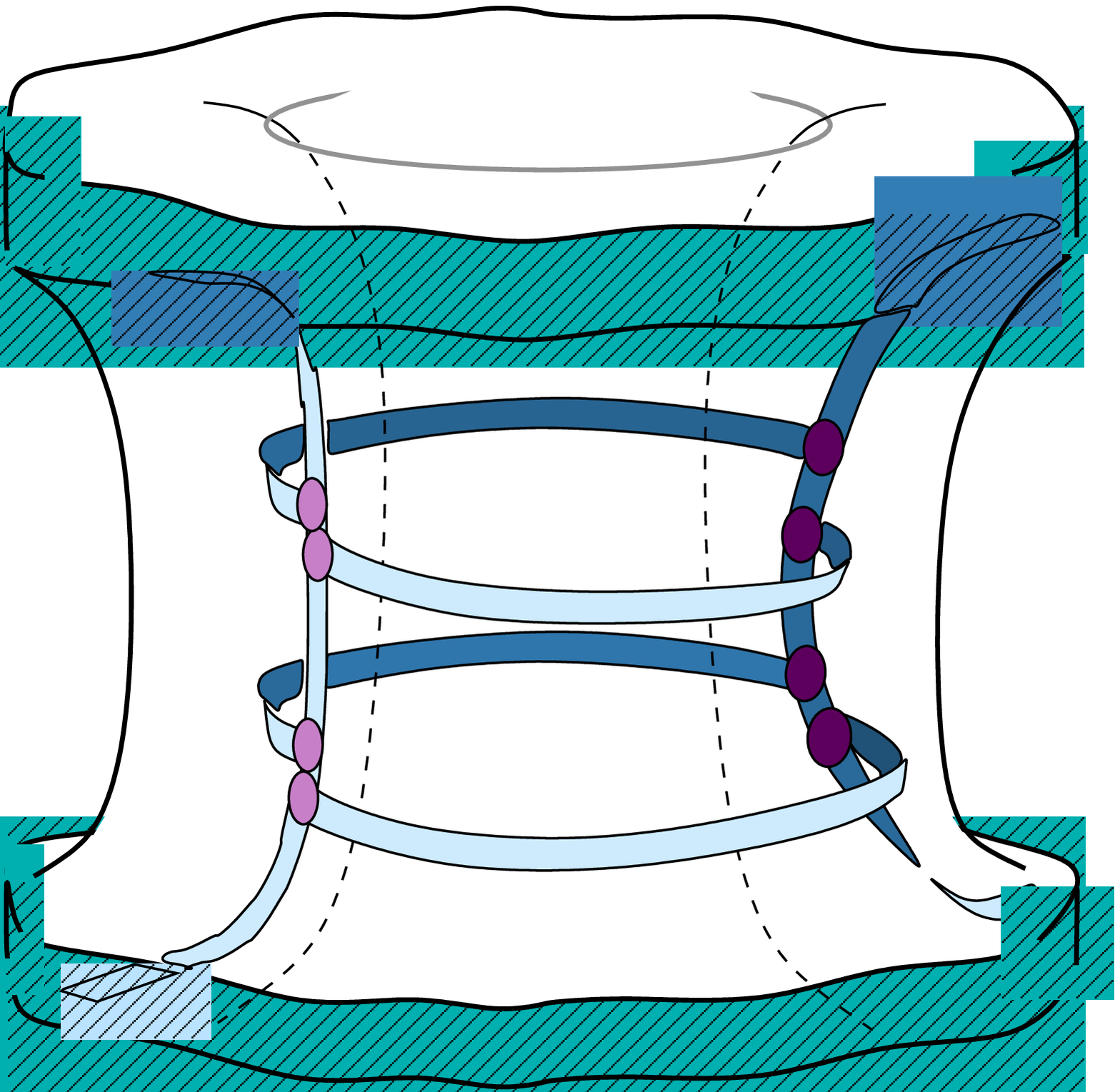}}}
       \setlength{\unitlength}{.3pt}
       \put(-40,335){\sse $M_1$}
       \put(-40, 43){\sse $M_2$}
       \put(75,75){\sse $A$}
       \put(204,198){\sse $A$}
       \put(334,249){\sse $A$}
       \put(188,107){\sse $A$}
       \setlength{\unitlength}{1pt}
       }
  \end{picture}
  \label{eq:app-bulk-1}
\ee
The two cylinders marked $M_1$ and $M_2$ indicate where the
part of $\M_A(\X,T)$ that is drawn connects to the rest of $\M_A(\X,T)$.
It is convenient to draw this fragment of cobordism
in `wedge representation' (see \cite[section 5.1]{fjfrs} for details),
\be
  \begin{picture}(288,208)(0,0)
  \put(20,100){(\ref{eq:app-bulk-1}) ~$=$ }
  \put(50,-1){
       \put(0,0){\scalebox{.3}{\includegraphics{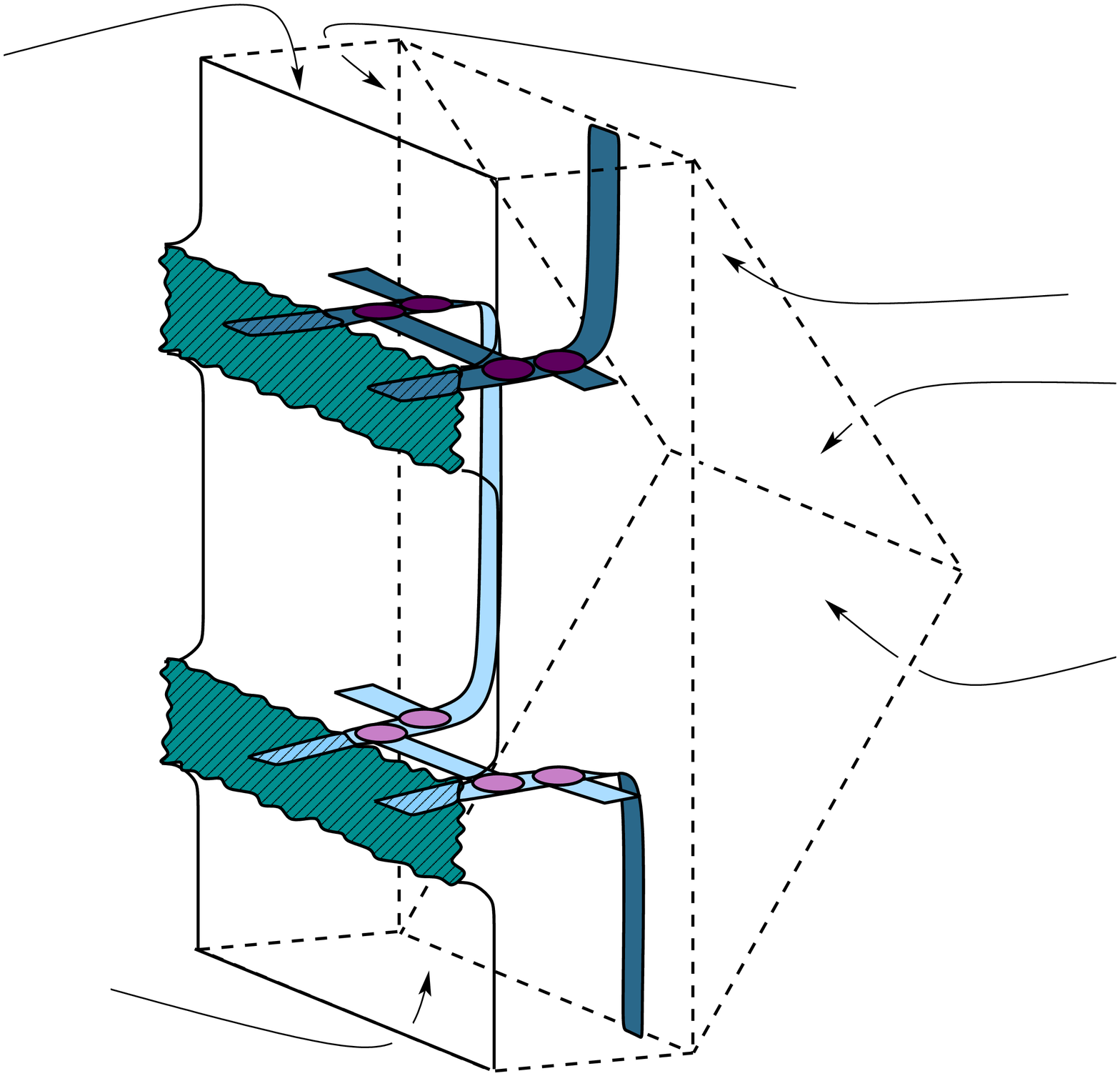}}}
       \setlength{\unitlength}{.3pt}
       \put(58,490){\sse $M_1$}
       \put(58,222){\sse $M_2$}
       \put(-26,650){\sse $B$}
       \put( 45, 52){\sse $B'$}
       \put(518,627){\sse $C$}
       \put(694,495){\sse $C'$}
       \put(725,440){\sse $D$}
       \put(725,264){\sse $D'$}
       \put(281,301){\sse $A$}
       \put(376, 83){\sse $A$}
       \put(400,502){\sse $A$}
       \setlength{\unitlength}{1pt}
       }
  \end{picture}
  \label{eq:app-bulk-2}
\ee
In this drawing, the surface labelled $B$ is to be identified
with $B'$ and similarly $C$ with $C'$ and $D$ with $D'$
It is then not too difficult
to see that (\ref{eq:app-bulk-2}) indeed describes the same
piece of cobordism as (\ref{eq:app-bulk-1}).

Next we note that the identity map on the vector space
$\tftc(T^2)$ assigned to a two-torus can be decomposed as
$$
  \begin{picture}(240,99)(0,0)
  \put(-8,52){$\id_{\tftc(T^2)}
  ~=~~ \displaystyle\sum_{i \in \II} \tftc\Big($}
  \put(227,55){$\Big)$}
  \put(100,0){
       \put(0,0){\scalebox{.25}{\includegraphics{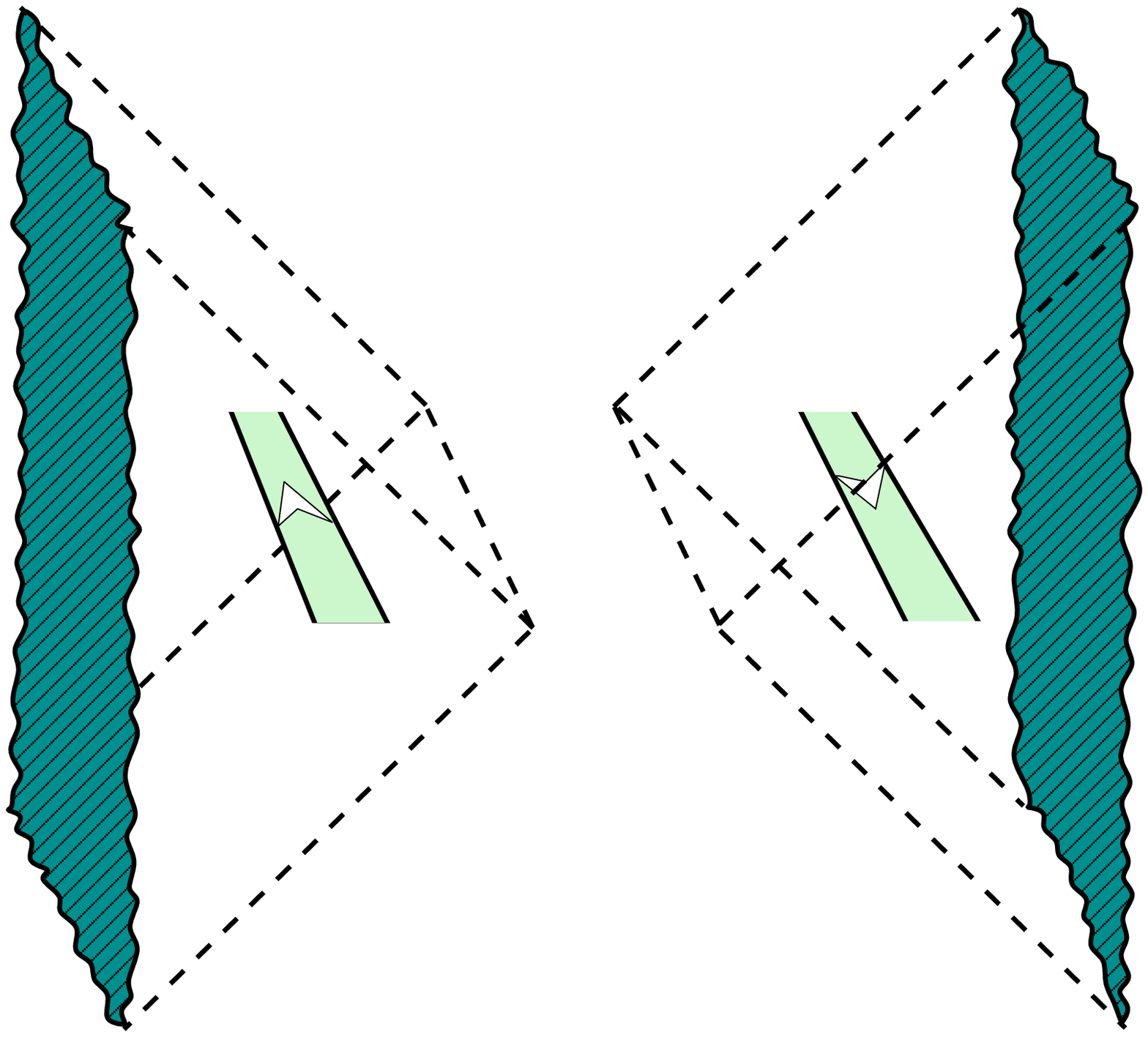}}}
       \setlength{\unitlength}{.25pt}
       \put(69,211){\sse $U_i$}
       \put(297,230){\sse $U_i$}
       \setlength{\unitlength}{1pt}
       }
  \end{picture}
$$
Applying this identity to (\ref{eq:app-bulk-2}) gives
(here we suppress the $\tftc(\,\cdots)$ that is to be taken for each cobordism)
\be
  \begin{picture}(200,262)(0,0)
  \put(-81,182){(\ref{eq:app-bulk-1}) $\,=\,
  \displaystyle\sum_{i \in \II}$}
  \put(-22,97){
       \put(-1,0){\scalebox{.24}{\includegraphics{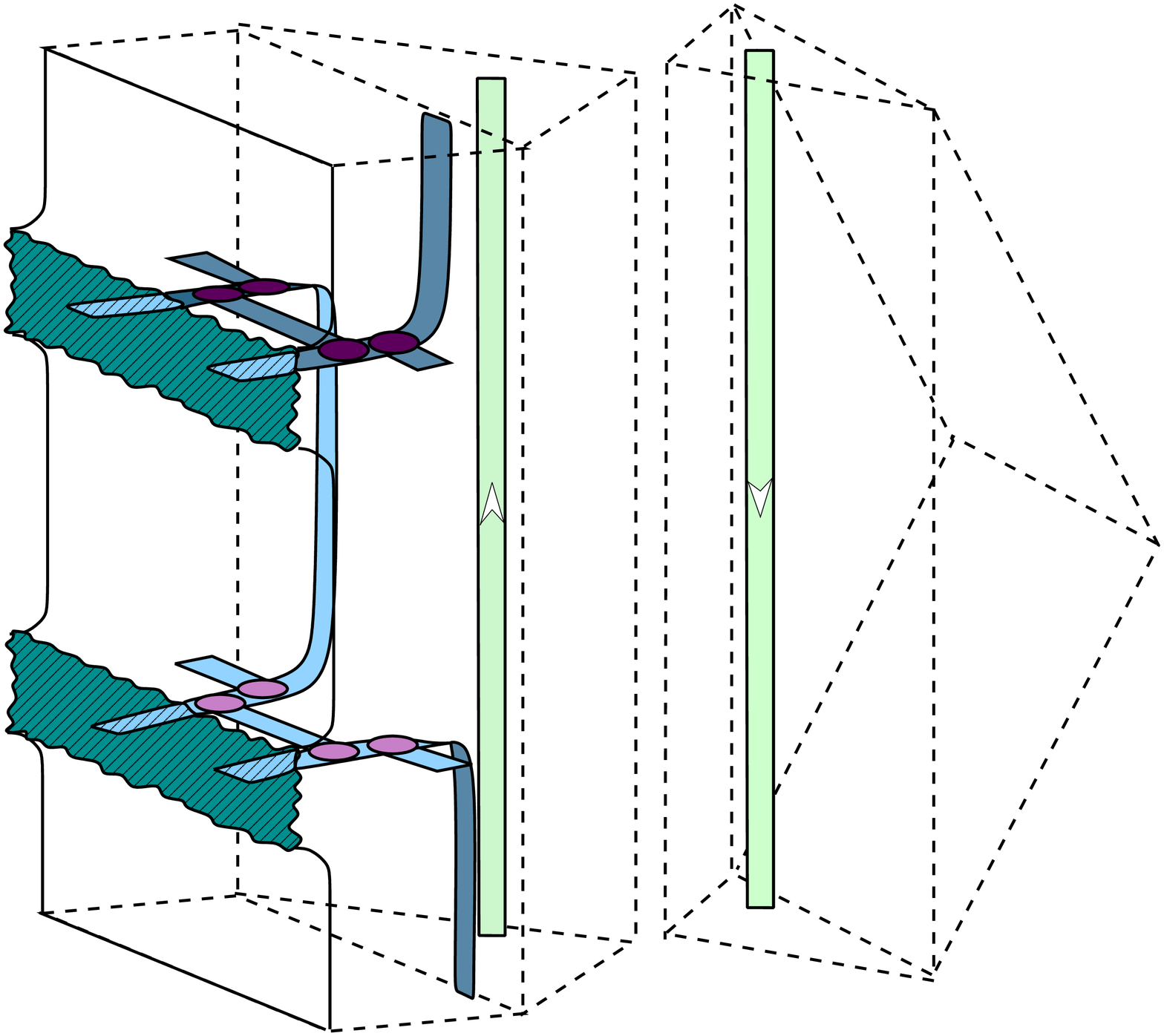}}}
       \setlength{\unitlength}{.25pt}
       \put(-47,472){\sse $M_1$}
       \put(-47,220){\sse $M_2$}
       \put(167,278){\sse $A$}
       \put(255,102){\sse $A$}
       \put(261,288){\sse $U_i$}
       \put(494,358){\sse $U_i$}
       \setlength{\unitlength}{1pt}
       }
  \put(132,88){$=\, \displaystyle\sum_{i\in\II} S_{0,0}^{}\,d_i^{}$}
  \put(188,0){
       \put(0,0){\scalebox{.2}{\includegraphics{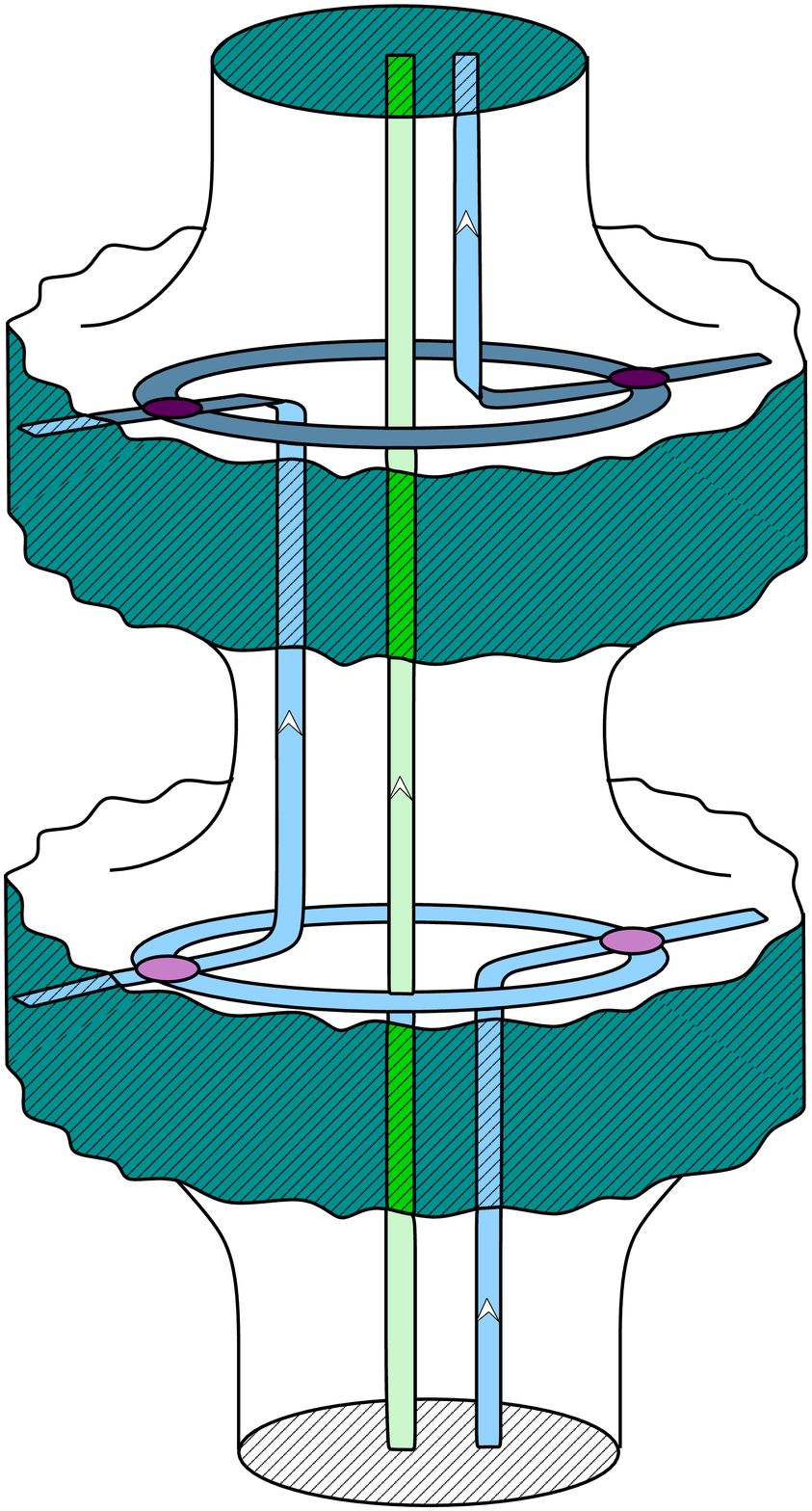}}}
       \setlength{\unitlength}{.2pt}
       \put(-58,642){\sse $M_1$}
       \put(-58,280){\sse $M_2$}
       \put(120,358){\sse $A$}
       \put(248,408){\sse $U_i$}
       \put(290,727){\sse $A$}
       \setlength{\unitlength}{1pt}
       }
  \end{picture}
  \label{eq:app-bulk-3}
\ee
with $d_i\,{:=}\,\dim(U_i)$. For the
second step, note that the connected component that contains only a
$U_i$-ribbon is in fact an $S^3$, and its invariant is $S_{0,0} \dim(U_i)$.

The picture (\ref{eq:app-bulk-1}) shows a fragment of the cobordism
for $\Cf(\X)$ close to the circle along which $\X$ is cut.
The right hand side of \eqref{eq:app-bulk-3} is nothing but the relevant
fragment of the cobordism obtained from applying ${\rm tr}_{\rm last}$ 
(as given in (\ref{eq:trlast}) with $\M_f$ as in (\ref{eq:Mf-trlast-bulk}))
to $\Cf({\rm cut}_f(\X))$. The relevant fragment of the cobordism for
$\Cf({\rm cut}_f(\X))$, in turn, is obtained by considering the second
picture in (\ref{eq:cut-bulk-double}) and inserting the ribbon graphs
(\ref{eq:in-out-bulk}) arising from an incoming and outgoing closed state 
boundary. The prefactor on the right hand side of \eqref{eq:app-bulk-3}
is produced by evaluating the two morphisms $\omega$ in 
(\ref{eq:in-out-bulk}). This establishes property {\bf(A2)}, i.e.\ that 
$\Cf(\X) \eq {\rm tr}_{\rm last}(\Cf({\rm cut}_f(\X) )$,
also for cutting along a circle.


\def\aff           {affine Lie algebra}
\def\alg           {algebra}
\def\app           {application}
\def\asis          {apparent singularities}
\def\auto          {automorphism}
\def\bc            {boundary condition}
\def\Bc            {Boundary condition}
\def\BGG           {Bern\-stein\hy Gel\-fand\hy Gel\-fand}
\def\bh            {black hole}
\def\Bh            {Black hole}
\def\bozn          {bosonization}
\def\bs            {boundary state}
\def\Bs            {Boundary state}
\def\Calg          {$C^*$-\al\-ge\-bra}
\def\cat           {category}
\def\cats          {categories}
\def\ccr           {canonical commutation relation}
\def\cft           {conformal field theory}
\def\Cft           {Conformal field theory}
\def\CFT           {Conformal Field Theory}
\def\cfts          {conformal field theories}
\def\Cfts          {Conformal field theories}
\def\cgc           {Clebsch\hy Gor\-dan coef\-fi\-cient}
\def\Class         {Classification }
\def\class         {classification}
\def\coeff         {coefficient}
\def\cocon         {coset construction}
\def\como          {coset model}
\def\compac        {compactification}
\def\complex       {{\dl C}}
\def\comrel        {commutation relation}
\def\Con           {Conformal }
\def\con           {conformal }
\def\conly         {conformally }
\def\Conly         {Conformally }
\def\corfu         {correlation function}
\def\Corfu         {Correlation function}
\def\cua           {current algebra}
\def\CY            {Cala\-bi\hy Yau }
\def\CYm           {Cala\-bi\hy Yau manifold}
\def\deq           {differential equation}
\def\dif           {differential}
\def\Dif           {Differential}
\def\dim           {dimension}
\def\dimred        {dimensional reduction}
\def\dimreg        {dimensional regularization}
\def\driso         {Drin\-feld\hy Soko\-lov }
\def\dyd           {Dynkin diagram}
\def\emt           {energy-momentum tensor}
\def\enva          {enveloping algebra}
\def\eq            {equa\-tion}
\def\findim        {fi\-ni\-te-di\-men\-si\-o\-nal}
\def\frob          {Fro\-be\-ni\-us algebra}
\def\fsi           {Fro\-be\-ni\-us\hy Schur indicator}
\def\ft            {field theory}
\def\fts           {field theories}
\def\furu          {fusion rule}
\def\Furu          {Fusion rule}
\def\gemo          {Gepner model}
\def\gkma          {generalized Kac\hy Moody algebra}
\def\gt            {gauge theory}
\def\gts           {gauge theories}
\def\gz            {generalized }
\def\Gz            {Generalized }
\def\Ham           {Ha\-mil\-to\-nian}
\def\Hamred        {Ha\-mil\-to\-nian reduction}
\def\hecke         {Hecke algebra}
\def\Hi            {Hilbert}
\def\hisc          {Hilbert scheme}
\def\hisp          {Hilbert space}
\def\hopf          {Hopf algebra}
\def\hw            {highest weight}
\def\Hw            {Highest weight}
\def\hy            {$\mbox{-\hspace{-.66 mm}-}$\linebreak[0]}
\newcommand\icdgmtp[1] {\rm {#1} International Conference on Differential
                   Geometric Methods in Theoretical Physics}
\def\ICM           {International Congress of Mathematicians }
\def\ide           {identification}
\def\Infdim        {Infinite-dimensional }
\def\infdim        {infinite-dimensional}
\def\intro         {introduction }
\def\Intro         {Introduction }
\def\inv           {invariance}
\def\irrep         {irreducible representation}
\def\KM            {Kac\hy Moody }
\def\kma           {Kac\hy Moody algebra}
\def\KN            {Krichever\hy Novikov }
\def\kna           {Krichever\hy Novikov algebra}
\def\Lag           {La\-gran\-gian}
\def\lcft          {logarithmic \cft}
\def\lcfts         {logarithmic \cfts}
\def\LG            {Lan\-dau\hy Ginz\-burg }
\def\lie           {Lie algebra}
\def\Lie           {Lie group}
\def\lmslns        {London Math.\ Soc.\ Lecture Note Series \# }
\def\lode          {linear ordinary differential equation}
\def\lrc           {Litt\-le\-wood\hy Ri\-chard\-son coefficient}
\def\maccor        {McKay correspondence}
\def\mimo          {minimal model}
\def\miss          {\query{\bf title miss}}
\def\Miss          {\query{\bf ref miss}}
\newcommand\misS[1]{\query{\bf #1}}
\def\Modinv        {Modular invarian}
\def\modinv        {modular invarian}
\def\mosp          {moduli space}
\def\mr            {reviewed in Math.\wb Reviews }
\def\MR            {Math.\wb Reviews}
\def\ncg           {noncommutative geometry}
\def\ncG           {non-com\-mu\-ta\-tive geometry}
\def\Ncg           {Noncommutative geometry}
\def\nlsm          {nonlinear sigma model}
\def\nlSm          {nonlinear $\sigma$ model}
\def\nLsm          {non-linear sigma model}
\def\nLsM          {non-linear sigma-model}
\def\nLSM          {non-linear $\sigma$-model}
\def\nn            {$N\,{=}\,2$ }
\def\NS            {Ne\-veu\hy Schwarz }
\def\oa            {operator algebra}
\def\onedim        {one-dimen\-sio\-nal}
\def\opa           {operator product algebra}
\def\ope           {operator product expansion}
\def\paf           {parafermion}
\def\parfu         {partition function}
\def\pf            {primary field}
\def\popo          {Poincar\'e polynomial}
\def\Q             {Quantum }
\def\q             {quantum }
\def\qdim          {quantum dimension}
\def\qft           {quantum field theory}
\def\Qft           {Quantum field theory}
\def\QFT           {Quantum Field Theory}
\def\qfts          {quantum field theories}
\def\qhe           {quantum Hall effect}
\def\qm            {quantum mechanics}
\def\qg            {quantum group}
\def\qgr           {quantum gravity}
\def\qsn           {quantisation}
\def\qzd           {quantized }
\def\qzn           {quantization}
\def\Qzn           {Quantization }
\def\rcfts         {ra\-ti\-o\-nal conformal field theory}
\def\reg           {regularization}
\def\ren           {renormaliz}
\def\Ren           {Renormaliz}
\def\reno          {renormalization }
\def\renO          {renormalisation }
\def\renog         {renormalization group}
\def\Renog         {Renormalization group }
\def\rep           {representation}
\def\Rep           {Representation}
\def\repth         {representation theory}
\def\Repth         {Representation theory}
\def\RI            {Riemann}
\def\RIsu          {Riemann surface}
\def\rmproc        {{\rm Proceedings of the }}
\def\sgts          {supersymmetric gauge theories}
\def\slie          {Lie su\-per\-al\-ge\-bra}
\def\ssi           {semisimple}
\def\ssI           {semi-simple}
\def\stc           {statistic}
\def\stt           {string theory}
\def\stts          {string theories}
\def\suco          {superconformal }
\def\sugra         {supergravity}
\def\suse          {superselection sector}
\def\susic         {supersymmetric }
\def\susiC         {supersymmetric}
\def\Susic         {Supersymmetric }
\def\susy          {supersymmetry}
\def\Susy          {Supersymmetry }
\def\SW            {Sei\-berg\hy Wit\-ten }
\def\sym           {symmetry}
\def\syms          {sym\-me\-tries}
\def\sYMt          {supersymmetric Yang\hy Mills theory}
\def\sYMts         {supersymmetric Yang\hy Mills theories}
\def\tc            {tensor category}
\def\tcs           {tensor categories}
\def\tft           {topological field theory}
\def\tfts          {topological field theories}
\def\Tft           {Topological field theory}
\def\Tfts          {Topological field theories}
\def\top           {topological }
\def\Top           {Topological }
\def\tqft          {topological quantum field theory}
\def\tqfts         {topological quantum field theories}
\def\trfo          {transformation}
\def\twodim        {two-di\-men\-sio\-nal}
\def\ttt           {To\-mi\-ta\hy Ta\-ke\-sa\-ki theory}
\def\uenva         {universal enveloping algebra}
\def\va            {Vi\-ra\-so\-ro algebra}
\def\vNa           {von Neumann algebra}
\def\voa           {vertex operator algebra}
\def\Voa           {Vertex operator algebra}
\def\vop           {vertex operator}
\def\Vop           {Vertex operator}
\def\YM            {Yang\hy Mills }
\def\YMt           {Yang\hy Mills theory}
\def\YMts          {Yang\hy Mills theories}
\def\WZ            {Wess\hy Zu\-mi\-no }
\def\WZW           {Wess\hy Zu\-mi\-no\hy Wit\-ten }
\def\wzwm          {WZW model}
\def\wzwt          {WZW theory}
\def\wzwts         {WZW the\-o\-ries}
\def\ybe           {Yang\hy Bax\-ter equation}
\def\zb            {reviewed in Zentralblatt No.\ }
\renewcommand\zb[1]{}
\newcommand\zB[2]{}
\renewcommand\mr[1]{}
\newcommand\mR[2]{}



\newcommand\wb{\,\linebreak[0]} \def\wB {$\,$\wb}
 \def\Bi              {\bibitem}
 \newcommand\JO[6]    {{\em #6}, {#1} {#2} ({#3}), {#4--#5} }
 \newcommand\J[7]     {{\em #7}, {#1} {#2} ({#3}), {#4--#5} {{\tt [#6]}}}
 \newcommand\K[6]     {{\em #6}, {#1} {#2} ({#3}), {#4} {{\tt [#5]}}}
 \newcommand\PhD[2]   {{\em #2\/}, Ph.D.\ thesis (#1)}
 \newcommand\Prep[2]  {{\em #2\/}, pre\-print {#1}}
 \newcommand\Pret[2]  {{\em #2}, pre\-print {\tt #1}}
 \newcommand\BOOK[4]  {{\sl #1\/} ({#2}, {#3} {#4})}
 \newcommand\inBO[7]  {{\em #7\/}, in:\ {\sl #1}, {#2}\ ({#3}, {#4} {#5}), p.\ {#6}}
 \newcommand\iNBO[7]  {{\em #7\/}, in:\ {\sl #1} ({#3}, {#4} {#5}) }
 \newcommand\Erra[3]  {\,[{\em ibid.}\ {#1} ({#2}) {#3}, {\em Erratum}]}
 \def\A     {Algebra }
 \def\dim   {dimension}
 \def\jf    {J.\ Fuchs}
 \def\adma  {Adv.\wb Math.}
 \def\anop  {Ann.\wb Phys.}
 \def\aspm  {Adv.\wb Stu\-dies\wB in\wB Pure\wB Math.}
 \def\atmp  {Adv.\wb Theor.\wb Math.\wb Phys.}
 \def\cocm  {Com\-mun.\wb Con\-temp.\wb Math.}
 \def\coia  {Com\-mun.\wB in\wB Algebra}
 \def\coma  {Con\-temp.\wb Math.}
 \def\comp  {Com\-mun.\wb Math.\wb Phys.}
 \def\cpma  {Com\-pos.\wb Math.}
 \def\duke  {Duke\wB Math.\wb J.}
 \def\fiic  {Fields\wB In\-sti\-tu\-te\wB Commun.}
 \def\fiiC  {Fields\wB Institute Commun.}
 \def\foph  {Fortschritte d.\wb Phys.}
 \def\gafa  {Geom.\wB and\wB Funct.\wb Anal.}
 \def\gatm  {Geom.\wB and\wB Topol.\wb Monogr.} 
 \def\hhaa  {Homol.\wb Homot.\wb Appl.}   
 \def\ijmb  {Int.\wb J.\wb Mod.\wb Phys.\ B}
 \def\ijmp  {Int.\wb J.\wb Mod.\wb Phys.\ A}
 \def\jams  {J.\wb Amer.\wb Math.\wb Soc.}
 \def\jgap  {J.\wb Geom.\wB and\wB Phys.}
 \def\jhep  {J.\wb High\wB Energy\wB Phys.}
 \def\jktr  {J.\wB Knot\wB Theory\wB and\wB its\wB Ramif.}
 \def\jlms  {J.\wB London\wB Math.\wb Soc.}
 \def\joal  {J.\wB Al\-ge\-bra}
 \def\jomp  {J.\wb Math.\wb Phys.}
 \def\jopa  {J.\wb Phys.\ A}
 \def\josp  {J.\wb Stat.\wb Phys.}
 \def\jpaa  {J.\wB Pure\wB Appl.\wb Alg.}
 \def\lemp  {Lett.\wb Math.\wb Phys.}
 \def\maan  {Math.\wb Annal.}
 \def\mams  {Memoirs\wB Amer.\wb Math.\wb Soc.}
 \def\mpla  {Mod.\wb Phys.\wb Lett.\ A}
 \def\nupb  {Nucl.\wb Phys.\ B}
 \def\pajm  {Pa\-cific\wB J.\wb Math.}
 \def\pams  {Proc.\wb Amer.\wb Math.\wb Soc.}
 \def\phep  {Proc.\wb HEP$\!$}
 \def\phla  {Phys.\wb Lett.\ A}
 \def\phlb  {Phys.\wb Lett.\ B}
 \def\phrd  {Phys.\wb Rev.\ D}
 \def\phre  {Phys.\wb Rev.\ E}
 \def\phrl  {Phys.\wb Rev.\wb Lett.}
 \def\phrp  {Phys.\wb Rep.}
 \def\pnas  {Proc.\wb Natl.\wb Acad.\wb Sci.\wb USA}
 \def\prtp  {Progr.\wb Theor.\wb Phys.}
 \def\rims  {Publ.\wB RIMS}
 \def\rpip  {Rep.\wb Prog.\wB in\wB Phys.}
 \def\rvmp  {Rev.\wb Math.\wb Phys.}
 \def\slnm  {Sprin\-ger\wB Lecture\wB Notes\wB in\wB Ma\-the\-matics }
 \newcommand\Slnm[1] {{\rm[\slnm\ #1]}}
 \def\taac  {Theory\wB and\wB \wB Appl.\wb Cat.}
 \def\tams  {Trans.\wb Amer.\wb Math.\wb Soc.}
 \def\thmp  {Theor.\wb Math.\wb Phys.}
 \def\topo  {Topology}
 \def\trgr  {Trans\-form. Groups}
 \def\AMS    {{American Mathematical Society}}
 \def\AP     {{Academic Press}}
 \def\CUP    {{Cambridge University Press}}
 \def\IPC    {{International Press Company}}
 \def\NH     {{North Holland Publishing Company}}
 \def\KLU    {{Kluwer Academic Publishers}}
 \def\PUP    {{Princeton University Press}}
 \def\PL     {{Plenum Press}}
 \def\SV     {{Sprin\-ger Ver\-lag}}
 \def\WI     {{Wiley Interscience}}
 \def\WS     {{World Scientific}}
 \def\Ad     {{Amsterdam}}
 \def\Be     {{Berlin}}
 \def\Ca     {{Cambridge}}
 \def\Do     {{Dordrecht}}
 \def\pR     {{Princeton}}
 \def\PR     {{Providence}}
 \def\Si     {{Singapore}}
 \def\NY     {{New York}}
   \newcommand\ncmp[2] {\inBO{IXth International Congress on
              Mathematical Physics} {B.\ Simon, A.\ Truman, and I.M.\ Davis,
              eds.} {Adam Hilger}{Bristol}{1989} {{#1}}{{#2}} } 
   \def\aste  {Ast\'e\-ris\-que}
   \def\BIR    {{Birk\-h\"au\-ser}}
   \def\Bo     {{Boston}}
   \def\ihes  {Publ.\wb Math.\wB I.H.E.S.}    

   \newcommand\bi\bibitem

   \def\inma  {Invent. math.}
   \def\taia  {Topology \wB and\wB its\wB Appl.}
   \def\joot  {J.\wB Operator\wB Theory}
   \def\jetl  {Sov.\wb Phys.\wB JETP\wB Lett.}
  \def\phrv  {Phys.\wb Rev.}
      \def\AW     {{Addi\-son\hy Wes\-ley}}
   \def\anma  {Ann.\wb Math.}
      \def\npbp  {Nucl.\wb Phys.\ B (Proc.\wb Suppl.)}
   \def\ptrs  {Phil.\wb Trans.\wb Roy.\wb Soc.\wB Lon\-don}

\bigskip
\bibliographystyle{amsalpha}

\end{document}